\documentclass[12pt,a4paper,leqno]{amsart}
\usepackage{amsmath}
\usepackage{amssymb}
\usepackage{bbm}
\usepackage{amsfonts}

\usepackage{amsthm}
\newtheorem{definition}{Definition}[section]
\newtheorem{theorem}[definition]{Theorem}
\newtheorem{proposition}[definition]{Proposition}
\newtheorem{remark}[definition]{Remark}
\newtheorem{corollary}[definition]{Corollary}
\newtheorem{lemma}[definition]{Lemma}

\usepackage{indentfirst}
\usepackage{multirow}
\numberwithin{equation}{section}
\allowdisplaybreaks[4] 

\usepackage{geometry}
\usepackage{tikz}

\usepackage{authblk} 
\usepackage[foot]{amsaddr} 

\title{Virial Theorem and Its Applications in Instability of Two-Phase Water-Wave}

\author{Haocheng Yang$^{1,2}$}
\thanks{$^{1}$Ecole Normale Supérieure Paris-Saclay, CNRS Centre Borelli UMR9010, 4 Avenue des Sciences, F-91190 Gif-sur-Yvette}
\thanks{$^{2}$Universit{\'e} Paris XIII (Sorbonne Paris-Nord), LAGA, CNRS (UMR 7539), 99 Avenue J.-B. Cl{\'e}ment, F-93430 Villetaneuse}

\usepackage[backend=biber, isbn=false, url=false, doi=false, maxnames=5]{biblatex}
\renewbibmacro{in:}{}
\usepackage{enumerate}
\usepackage{verbatim} 

\usepackage[colorlinks,linkcolor=blue,citecolor=red,urlcolor=black]{hyperref}

\begin{document}
	
	\pagenumbering{arabic}
	
	\newcommand{\Op}[2][]{\operatorname{Op}^{#1}\left(#2\right)}
	\newcommand{\Supp}[1]{\operatorname{Supp} #1}
	\newcommand{\Real}{\operatorname{Re}}
	\newcommand{\Imag}{\operatorname{Im}}
	\newcommand{\Id}{\mathrm{Id}}
	\newcommand{\B}{\mathcal{B}}
	\newcommand{\V}{\mathcal{V}}
	\newcommand{\diver}{\operatorname{div}}
	\newcommand{\curl}{\operatorname{curl}}
	\newcommand{\BBra}[1]{[\![ #1 ]\!]}
	
	\newcommand{\R}{\mathbb{R}}
	\newcommand{\C}{\mathbb{C}}
	\newcommand{\N}{\mathbb{N}}
	\newcommand{\T}{\mathbb{T}}
	\newcommand{\Z}{\mathbb{Z}}
	
	\renewcommand{\le}{\leqslant}
	\renewcommand{\ge}{\geqslant}
	
	\newcommand{\mar}[1]{\marginpar{\textcolor{red}{#1}}} 

	\begin{abstract}
		In this paper, we analyze the dynamics of two layers of immiscible, inviscid, and incompressible potential fluids through a fully nonlinear system. Our goal is to establish a virial theorem and prove the algebric growth of slope and curvature of the interface over time when the fluid below is no denser than the one above. These phenomena, known as Rayleigh-Taylor instability, will be proved for a broad class of regular initial data, including the case of 2D overlapping interface.
	\end{abstract}

	\maketitle
	
	\tableofcontents

	\section{Introduction}\label{Sect:intro}
	
	\subsection{Two-phase water-wave}\label{subsect-intro:WW}
	
	Consider two fluids occupying the domains $\Omega^\pm(t)\subset \T^d\times]-H^+,H^-[$ separated by a hypersurface $\Sigma(t)$ without self-intersection, parametrized by
	\begin{equation}\label{eq-intro-ww:DefFrBdy}
		\gamma(t)\colon \T^d \to \T^d\times]-H^+,H^-[,\quad \hat{s}\mapsto \gamma(t,\hat{s}).
	\end{equation}
	Here $H^\pm\in]0,+\infty]$ stands for the depth and we denote by 
	\begin{equation}\label{eq-intro-ww:DefSolBdy}
		\Gamma^\pm=\{(x,y)\in\T^d\times\R:y=\mp H^\pm\}
	\end{equation}
	the solid boundaries (or the ``bottoms''). When the depth $H^+=+\infty$ or $H^-=+\infty$, we do not consider the corresponding solid boundary. Note that, in this article, the dimension $d$ is chosen to be $1$ or $2$ for physical relevance, while most of the arguments hold true for general dimensions.
	\begin{figure}[h]\label{fig-intro-ww:NonGraph}
		\centering
		\begin{tikzpicture}[x=0.6pt,y=0.6pt,yscale=-1,xscale=1]
			
			\draw  (42,182.67) -- (467.73,182.67)(71.73,7.8) -- (71.73,299.8) (460.73,177.67) -- (467.73,182.67) -- (460.73,187.67) (66.73,14.8) -- (71.73,7.8) -- (76.73,14.8)  ;
			\draw    (71.2,54.8) -- (431.2,54.8) ;
			\draw    (71.2,264.8) -- (431.2,264.8) ;
			\draw    (71.2,44.8) -- (81.2,54.8) ;
			\draw    (81.34,44.8) -- (91.34,54.8) ;
			\draw    (91.2,44.8) -- (101.2,54.8) ;
			\draw    (101.77,44.8) -- (111.77,54.8) ;
			\draw    (111.2,44.8) -- (121.2,54.8) ;
			\draw    (121.34,44.8) -- (131.34,54.8) ;
			\draw    (131.2,44.8) -- (141.2,54.8) ;
			\draw    (141.77,44.8) -- (151.77,54.8) ;
			\draw    (151.8,44.8) -- (161.8,54.8) ;
			\draw    (161.94,44.8) -- (171.94,54.8) ;
			\draw    (171.8,44.8) -- (181.8,54.8) ;
			\draw    (182.37,44.8) -- (192.37,54.8) ;
			\draw    (191.8,44.8) -- (201.8,54.8) ;
			\draw    (201.94,44.8) -- (211.94,54.8) ;
			\draw    (211.8,44.8) -- (221.8,54.8) ;
			\draw    (222.37,44.8) -- (232.37,54.8) ;
			\draw    (231.56,44.8) -- (241.56,54.8) ;
			\draw    (241.71,44.8) -- (251.71,54.8) ;
			\draw    (251.56,44.8) -- (261.56,54.8) ;
			\draw    (262.14,44.8) -- (272.14,54.8) ;
			\draw    (271.56,44.8) -- (281.56,54.8) ;
			\draw    (281.71,44.8) -- (291.71,54.8) ;
			\draw    (291.56,44.8) -- (301.56,54.8) ;
			\draw    (302.14,44.8) -- (312.14,54.8) ;
			\draw    (312.16,44.8) -- (322.16,54.8) ;
			\draw    (322.31,44.8) -- (332.31,54.8) ;
			\draw    (332.16,44.8) -- (342.16,54.8) ;
			\draw    (342.74,44.8) -- (352.74,54.8) ;
			\draw    (352.16,44.8) -- (362.16,54.8) ;
			\draw    (362.31,44.8) -- (372.31,54.8) ;
			\draw    (372.16,44.8) -- (382.16,54.8) ;
			\draw    (382.74,44.8) -- (392.74,54.8) ;
			\draw    (391.2,44.8) -- (401.2,54.8) ;
			\draw    (401.34,44.8) -- (411.34,54.8) ;
			\draw    (411.2,44.8) -- (421.2,54.8) ;
			\draw    (421.77,44.8) -- (431.77,54.8) ;
			\draw    (71.2,264.8) -- (81.2,274.8) ;
			\draw    (81.34,264.8) -- (91.34,274.8) ;
			\draw    (91.2,264.8) -- (101.2,274.8) ;
			\draw    (101.77,264.8) -- (111.77,274.8) ;
			\draw    (111.2,264.8) -- (121.2,274.8) ;
			\draw    (121.34,264.8) -- (131.34,274.8) ;
			\draw    (131.2,264.8) -- (141.2,274.8) ;
			\draw    (141.77,264.8) -- (151.77,274.8) ;
			\draw    (151.8,264.8) -- (161.8,274.8) ;
			\draw    (161.94,264.8) -- (171.94,274.8) ;
			\draw    (171.8,264.8) -- (181.8,274.8) ;
			\draw    (182.37,264.8) -- (192.37,274.8) ;
			\draw    (191.8,264.8) -- (201.8,274.8) ;
			\draw    (201.94,264.8) -- (211.94,274.8) ;
			\draw    (211.8,264.8) -- (221.8,274.8) ;
			\draw    (222.37,264.8) -- (232.37,274.8) ;
			\draw    (231.56,264.8) -- (241.56,274.8) ;
			\draw    (241.71,264.8) -- (251.71,274.8) ;
			\draw    (251.56,264.8) -- (261.56,274.8) ;
			\draw    (262.14,264.8) -- (272.14,274.8) ;
			\draw    (271.56,264.8) -- (281.56,274.8) ;
			\draw    (281.71,264.8) -- (291.71,274.8) ;
			\draw    (291.56,264.8) -- (301.56,274.8) ;
			\draw    (302.14,264.8) -- (312.14,274.8) ;
			\draw    (312.16,264.8) -- (322.16,274.8) ;
			\draw    (322.31,264.8) -- (332.31,274.8) ;
			\draw    (332.16,264.8) -- (342.16,274.8) ;
			\draw    (342.74,264.8) -- (352.74,274.8) ;
			\draw    (352.16,264.8) -- (362.16,274.8) ;
			\draw    (362.31,264.8) -- (372.31,274.8) ;
			\draw    (372.16,264.8) -- (382.16,274.8) ;
			\draw    (382.74,264.8) -- (392.74,274.8) ;
			\draw    (391.2,264.8) -- (401.2,274.8) ;
			\draw    (401.34,264.8) -- (411.34,274.8) ;
			\draw    (411.2,264.8) -- (421.2,274.8) ;
			\draw    (421.77,264.8) -- (431.77,274.8) ;
			\draw [color={rgb, 255:red, 65; green, 117; blue, 5 }  ,draw opacity=1 ][line width=1.5]    (71.73,182.67) .. controls (96.09,180.95) and (142.6,171.55) .. (170.1,173.05) .. controls (197.6,174.55) and (215.6,205.55) .. (256.1,210.05) .. controls (296.6,214.55) and (326.1,202.55) .. (340.1,181.05) .. controls (354.1,159.55) and (346.2,148.35) .. (310.6,145.55) .. controls (275,142.75) and (264.85,129.55) .. (283.6,114.55) .. controls (302.35,99.55) and (360.6,119.05) .. (385.6,144.55) .. controls (410.6,170.05) and (409.1,182.55) .. (429.28,182.96) ;
			\draw [color={rgb, 255:red, 65; green, 117; blue, 5 }  ,draw opacity=1 ][fill={rgb, 255:red, 65; green, 117; blue, 5 }  ,fill opacity=1 ]   (275.1,123.05) -- (238.35,102.53) ;
			\draw [shift={(236.6,101.55)}, rotate = 29.18] [color={rgb, 255:red, 65; green, 117; blue, 5 }  ,draw opacity=1 ][line width=0.75]    (10.93,-3.29) .. controls (6.95,-1.4) and (3.31,-0.3) .. (0,0) .. controls (3.31,0.3) and (6.95,1.4) .. (10.93,3.29)   ;
			\draw [color={rgb, 255:red, 74; green, 144; blue, 226 }  ,draw opacity=1 ]   (350.1,54.55) -- (350.38,23.05) ;
			\draw [shift={(350.4,21.05)}, rotate = 90.51] [color={rgb, 255:red, 74; green, 144; blue, 226 }  ,draw opacity=1 ][line width=0.75]    (10.93,-3.29) .. controls (6.95,-1.4) and (3.31,-0.3) .. (0,0) .. controls (3.31,0.3) and (6.95,1.4) .. (10.93,3.29)   ;
			\draw [color={rgb, 255:red, 245; green, 166; blue, 35 }  ,draw opacity=1 ]   (352.16,264.8) -- (352.39,298.05) ;
			\draw [shift={(352.4,300.05)}, rotate = 269.62] [color={rgb, 255:red, 245; green, 166; blue, 35 }  ,draw opacity=1 ][line width=0.75]    (10.93,-3.29) .. controls (6.95,-1.4) and (3.31,-0.3) .. (0,0) .. controls (3.31,0.3) and (6.95,1.4) .. (10.93,3.29)   ;
			
			\draw (79.2,1.87) node [anchor=north west][inner sep=0.75pt]    {$y$};
			\draw (415.93,189.37) node [anchor=north west][inner sep=0.75pt]    {$x\in \mathbb{T}{^{d}}$};
			\draw (28.4,44) node [anchor=north west][inner sep=0.75pt]    {$H^{-}$};
			\draw (20.3,257.3) node [anchor=north west][inner sep=0.75pt]    {$-H^{+}$};
			\draw (108.6,145.4) node [anchor=north west][inner sep=0.75pt]    {$\textcolor[rgb]{0.25,0.46,0.02}{\Sigma }\textcolor[rgb]{0.25,0.46,0.02}{(}\textcolor[rgb]{0.25,0.46,0.02}{t}\textcolor[rgb]{0.25,0.46,0.02}{)}\textcolor[rgb]{0.25,0.46,0.02}{=}\textcolor[rgb]{0.25,0.46,0.02}{\{}\textcolor[rgb]{0.25,0.46,0.02}{\gamma }\textcolor[rgb]{0.25,0.46,0.02}{(}\textcolor[rgb]{0.25,0.46,0.02}{t,\cdot }\textcolor[rgb]{0.25,0.46,0.02}{)}\textcolor[rgb]{0.25,0.46,0.02}{\}}$};
			\draw (142.3,84.2) node [anchor=north west][inner sep=0.75pt]  [color={rgb, 255:red, 74; green, 144; blue, 226 }  ,opacity=1 ]  {$\Omega ^{-}( t)$};
			\draw (140.8,219.7) node [anchor=north west][inner sep=0.75pt]  [color={rgb, 255:red, 245; green, 166; blue, 35 }  ,opacity=1 ]  {$\Omega ^{+}( t)$};
			\draw (403.2,58.2) node [anchor=north west][inner sep=0.75pt]  [color={rgb, 255:red, 74; green, 144; blue, 226 }  ,opacity=1 ]  {$\Gamma ^{-}$};
			\draw (405.3,240.2) node [anchor=north west][inner sep=0.75pt]  [color={rgb, 255:red, 245; green, 166; blue, 35 }  ,opacity=1 ]  {$\Gamma ^{+}$};
			\draw (255.6,86.9) node [anchor=north west][inner sep=0.75pt]  [color={rgb, 255:red, 65; green, 117; blue, 5 }  ,opacity=1 ]  {$n$};
			\draw (364.1,11.4) node [anchor=north west][inner sep=0.75pt]  [color={rgb, 255:red, 74; green, 144; blue, 226 }  ,opacity=1 ]  {$n_{b}^{-}$};
			\draw (368.6,287.4) node [anchor=north west][inner sep=0.75pt]  [color={rgb, 255:red, 245; green, 166; blue, 35 }  ,opacity=1 ]  {$n_{b}^{+}$};

		\end{tikzpicture}
		
		\caption{Two phase water-wave}
	\end{figure}
	
	Since the fluids are assumed to be potential, inviscid, and incompressible, their velocity field $u^\pm$ is expressed as $u^\pm = \nabla_{x,y}\phi^\pm$ for some scalar potential $\phi^\pm$ and satisfies the Euler equation with divergence-free condition,
	\begin{equation}\label{eq:Euler}
		\left\{\begin{array}{ll}
			\rho^\pm\left(\partial_t u^\pm + u^\pm\cdot\nabla_{x,y} u^\pm\right) + \nabla_{x,y} \left(P^\pm + \rho^\pm gy\right) = 0, & \text{in }\Omega^\pm(t), \\ [0.5ex]
			\diver_{x,y}u^\pm = 0, & \text{in }\Omega^\pm(t), \\[0.5ex]
			u^\pm=\nabla_{x,y}\phi^\pm, & \text{in }\Omega^\pm(t),
		\end{array}\right.
	\end{equation}
	where $\rho^\pm\ge0$ (with $\rho^++\rho^->0$ to avoid trivial case) is constant density, $P^\pm$ is the pressure, and $g\ge0$ is the gravity acceleration. In order to complete the description of the system, we need to determine the boundary conditions on $\Sigma(t)$ and $\Gamma^\pm$. Since the fluids cannot go through the outer boundary $\{y=\mp H^\pm\}$, one has
	\begin{equation}\label{eq:bdy-cond-bottom}
		u^\pm \cdot n^\pm_b = 0,\ \ \text{on }\Gamma^\pm,
	\end{equation}
	where $n^\pm_b = (0,\mp 1)$ is the unit outward normal vector of the solid boundary $\Gamma^\pm$. Remark that, in infinite depth case $H^{+}=+\infty$ or $H^{-}=+\infty$, this condition should be replaced by
	\begin{equation}\label{eq:bdy-cond-bottom-infty-depth}
		\lim_{y\rightarrow-\infty} u^+=0,\ \ \text{or }\lim_{y\rightarrow+\infty} u^-=0,
	\end{equation}
	respectively. As for the interface $\Sigma(t)$, we assume that the normal component of the velocity field is continuous to avoid vacuum area, which is known as \textit{kinematic boundary condition},
	\begin{equation}\label{eq:bdy-cond-interface}
		\partial_t\gamma \cdot n = u^+\cdot n = u^-\cdot n,\ \ \text{on }\Sigma(t),
	\end{equation}
	where $n$ is the unit normal vector of interface. The last boundary condition is about the pressure $P^\pm$,
	\begin{equation}\label{eq:bdy-cond-pressure}
		P^+ - P^- = \sigma\kappa,\ \ \text{on }\Sigma(t),
	\end{equation}
	where $\sigma\ge0$ is the surface tension constant ($\sigma=0$ represents the case without surface tension), and $\kappa$ is the mean curvature.
	
	Since the velocity field $u^\pm$ can be expressed as the gradient $\nabla_{x,y}\phi^\pm$ of scalar potential $\phi^\pm$, the divergence-free condition implies
	\begin{equation}\label{eq-intro-ww:DefSclPot}
		\Delta_{x,y}\phi^\pm = 0,\ \ \text{in }\Omega^\pm(t).
	\end{equation}
	Consequently, the Euler equation \eqref{eq:Euler} yields the \textit{Bernoulli's equation},
	\begin{equation}\label{eq-intro-ww:Bernoulli}
		\rho^\pm\left(\partial_t \phi^\pm + \frac{|\nabla_{x,y}\phi^\pm|^2}{2}\right) + P^\pm + \rho^\pm gy = 0,\ \  \text{in }\Omega^\pm(t).
	\end{equation}
	Note that the right-hand side is generally a constant depending only on time and could be taken as zero by absorbing the constant in the definition of the pressure $P^\pm$. Now, we can state the main system of $(\gamma,\phi^\pm,P^\pm)$ to be studied throughout this article,
	\begin{equation}\label{eq-intro-ww:Main}
		\left\{\begin{array}{ll}
			\rho^\pm\left(\partial_t \phi^\pm + \frac12|\nabla_{x,y}\phi^\pm|^2\right) + P^\pm + \rho^\pm gy = 0, &\text{in }\Omega^\pm(t), \\[0.5ex]
			\Delta_{x,y}\phi^\pm = 0, &\text{in }\Omega^\pm(t), \\[0.5ex]
			n\cdot \partial_t\gamma = n \cdot \nabla_{x,y}\phi^+ = n \cdot \nabla_{x,y}\phi^-, &\text{on }\Sigma(t), \\[0.5ex]
			P^+ - P^- = \sigma\kappa, &\text{on }\Sigma(t), \\[0.5ex]
			\text{(finite depth) }n^\pm_b \cdot \nabla_{x,y}\phi^\pm = 0, &\text{on }\Gamma^\pm, \\
			\text{(infinite depth) }\lim_{y\to\mp\infty}\nabla_{x,y}\phi^\pm=0. &
		\end{array}\right.
	\end{equation}
	
	In fact, it is possible to determine the pressure $P^\pm$ from $(\gamma, \phi^\pm)$, via \textit{Craig-Sulem-Zakharov formulation}, cf. \cite{zakharov1968stability,craig1993numerical,benjamin1997reappraisal}. Roughly speaking, we compute the difference of $+$ and $-$ parts of the first equation in \eqref{eq-intro-ww:Main} and restrict this difference at the free boundary $\Sigma(t)$. Then the fourth equation in \eqref{eq-intro-ww:Main} could cancel the pressures. In this way, we obtain a self-contained system involving only $\gamma$ and the traces of $\phi^\pm$, from which we can restore $\phi^\pm$ by solving the Poisson equation $\Delta_{x,y}\phi^\pm=0$ with proper boundary conditions. We refer to \cite{craig2017hamiltonian} for the detailed computations. To avoid technical details, this formulation will not be used in this article.

	There are two special cases to be studied in this article. One is the \textit{graph case} where $d\ge 1$ and $\Sigma(t)$ is a graph of some function $\eta(t)\colon\T^d\to]-H^+,H^-[$. The parametrization $\gamma(t)$ can be fixed as $\gamma(t,x)=(x,\eta(t,x))$. 
	
	\begin{figure}[h]
		\centering
		
		\begin{tikzpicture}[x=0.65pt,y=0.65pt,yscale=-1,xscale=1]
			
			\draw  (21.27,126.33) -- (302.19,126.33)(40.89,12.25) -- (40.89,202.74) (295.19,121.33) -- (302.19,126.33) -- (295.19,131.33) (35.89,19.25) -- (40.89,12.25) -- (45.89,19.25)  ;
			\draw    (40.54,42.91) -- (278.09,42.91) ;
			\draw    (40.54,179.91) -- (278.09,179.91) ;
			\draw    (40.54,36.38) -- (47.14,42.91) ;
			\draw    (47.23,36.38) -- (53.83,42.91) ;
			\draw    (53.74,36.38) -- (60.33,42.91) ;
			\draw    (60.71,36.38) -- (67.31,42.91) ;
			\draw    (66.93,36.38) -- (73.53,42.91) ;
			\draw    (73.63,36.38) -- (80.22,42.91) ;
			\draw    (80.13,36.38) -- (86.73,42.91) ;
			\draw    (87.11,36.38) -- (93.7,42.91) ;
			\draw    (93.72,36.38) -- (100.32,42.91) ;
			\draw    (100.42,36.38) -- (107.01,42.91) ;
			\draw    (106.92,36.38) -- (113.52,42.91) ;
			\draw    (113.9,36.38) -- (120.49,42.91) ;
			\draw    (120.12,36.38) -- (126.72,42.91) ;
			\draw    (126.81,36.38) -- (133.41,42.91) ;
			\draw    (133.31,36.38) -- (139.91,42.91) ;
			\draw    (140.29,36.38) -- (146.89,42.91) ;
			\draw    (146.36,36.38) -- (152.95,42.91) ;
			\draw    (153.05,36.38) -- (159.65,42.91) ;
			\draw    (159.55,36.38) -- (166.15,42.91) ;
			\draw    (166.53,36.38) -- (173.13,42.91) ;
			\draw    (172.75,36.38) -- (179.35,42.91) ;
			\draw    (179.44,36.38) -- (186.04,42.91) ;
			\draw    (185.95,36.38) -- (192.55,42.91) ;
			\draw    (192.92,36.38) -- (199.52,42.91) ;
			\draw    (199.54,36.38) -- (206.14,42.91) ;
			\draw    (206.23,36.38) -- (212.83,42.91) ;
			\draw    (212.74,36.38) -- (219.34,42.91) ;
			\draw    (219.71,36.38) -- (226.31,42.91) ;
			\draw    (225.93,36.38) -- (232.53,42.91) ;
			\draw    (232.63,36.38) -- (239.23,42.91) ;
			\draw    (239.13,36.38) -- (245.73,42.91) ;
			\draw    (246.11,36.38) -- (252.7,42.91) ;
			\draw    (251.69,36.38) -- (258.29,42.91) ;
			\draw    (258.38,36.38) -- (264.98,42.91) ;
			\draw    (264.89,36.38) -- (271.49,42.91) ;
			\draw    (271.86,36.38) -- (278.46,42.91) ;
			\draw    (40.54,179.91) -- (47.14,186.43) ;
			\draw    (47.23,179.91) -- (53.83,186.43) ;
			\draw    (53.74,179.91) -- (60.33,186.43) ;
			\draw    (60.71,179.91) -- (67.31,186.43) ;
			\draw    (66.93,179.91) -- (73.53,186.43) ;
			\draw    (73.63,179.91) -- (80.22,186.43) ;
			\draw    (80.13,179.91) -- (86.73,186.43) ;
			\draw    (87.11,179.91) -- (93.7,186.43) ;
			\draw    (93.72,179.91) -- (100.32,186.43) ;
			\draw    (100.42,179.91) -- (107.01,186.43) ;
			\draw    (106.92,179.91) -- (113.52,186.43) ;
			\draw    (113.9,179.91) -- (120.49,186.43) ;
			\draw    (120.12,179.91) -- (126.72,186.43) ;
			\draw    (126.81,179.91) -- (133.41,186.43) ;
			\draw    (133.31,179.91) -- (139.91,186.43) ;
			\draw    (140.29,179.91) -- (146.89,186.43) ;
			\draw    (146.36,179.91) -- (152.95,186.43) ;
			\draw    (153.05,179.91) -- (159.65,186.43) ;
			\draw    (159.55,179.91) -- (166.15,186.43) ;
			\draw    (166.53,179.91) -- (173.13,186.43) ;
			\draw    (172.75,179.91) -- (179.35,186.43) ;
			\draw    (179.44,179.91) -- (186.04,186.43) ;
			\draw    (185.95,179.91) -- (192.55,186.43) ;
			\draw    (192.92,179.91) -- (199.52,186.43) ;
			\draw    (199.54,179.91) -- (206.14,186.43) ;
			\draw    (206.23,179.91) -- (212.83,186.43) ;
			\draw    (212.74,179.91) -- (219.34,186.43) ;
			\draw    (219.71,179.91) -- (226.31,186.43) ;
			\draw    (225.93,179.91) -- (232.53,186.43) ;
			\draw    (232.63,179.91) -- (239.23,186.43) ;
			\draw    (239.13,179.91) -- (245.73,186.43) ;
			\draw    (246.11,179.91) -- (252.7,186.43) ;
			\draw    (251.69,179.91) -- (258.29,186.43) ;
			\draw    (258.38,179.91) -- (264.98,186.43) ;
			\draw    (264.89,179.91) -- (271.49,186.43) ;
			\draw    (271.86,179.91) -- (278.46,186.43) ;
			\draw [color={rgb, 255:red, 65; green, 117; blue, 5 }  ,draw opacity=1 ][line width=1.5]    (40.89,126.33) .. controls (56.97,125.2) and (102.37,93.96) .. (122.95,102.05) .. controls (143.54,110.14) and (164.53,146.41) .. (199.5,146.08) .. controls (234.47,145.76) and (262.94,125.24) .. (277.46,126.54) ;
			\draw  (346.27,126.33) -- (627.19,126.33)(365.89,12.25) -- (365.89,202.74) (620.19,121.33) -- (627.19,126.33) -- (620.19,131.33) (360.89,19.25) -- (365.89,12.25) -- (370.89,19.25)  ;
			\draw    (365.54,42.91) -- (603.09,42.91) ;
			\draw    (365.54,179.91) -- (603.09,179.91) ;
			\draw    (365.54,36.38) -- (372.14,42.91) ;
			\draw    (372.23,36.38) -- (378.83,42.91) ;
			\draw    (378.74,36.38) -- (385.33,42.91) ;
			\draw    (385.71,36.38) -- (392.31,42.91) ;
			\draw    (391.93,36.38) -- (398.53,42.91) ;
			\draw    (398.63,36.38) -- (405.22,42.91) ;
			\draw    (405.13,36.38) -- (411.73,42.91) ;
			\draw    (412.11,36.38) -- (418.7,42.91) ;
			\draw    (418.72,36.38) -- (425.32,42.91) ;
			\draw    (425.42,36.38) -- (432.01,42.91) ;
			\draw    (431.92,36.38) -- (438.52,42.91) ;
			\draw    (438.9,36.38) -- (445.49,42.91) ;
			\draw    (445.12,36.38) -- (451.72,42.91) ;
			\draw    (451.81,36.38) -- (458.41,42.91) ;
			\draw    (458.31,36.38) -- (464.91,42.91) ;
			\draw    (465.29,36.38) -- (471.89,42.91) ;
			\draw    (471.36,36.38) -- (477.95,42.91) ;
			\draw    (478.05,36.38) -- (484.65,42.91) ;
			\draw    (484.55,36.38) -- (491.15,42.91) ;
			\draw    (491.53,36.38) -- (498.13,42.91) ;
			\draw    (497.75,36.38) -- (504.35,42.91) ;
			\draw    (504.44,36.38) -- (511.04,42.91) ;
			\draw    (510.95,36.38) -- (517.55,42.91) ;
			\draw    (517.92,36.38) -- (524.52,42.91) ;
			\draw    (524.54,36.38) -- (531.14,42.91) ;
			\draw    (531.23,36.38) -- (537.83,42.91) ;
			\draw    (537.74,36.38) -- (544.34,42.91) ;
			\draw    (544.71,36.38) -- (551.31,42.91) ;
			\draw    (550.93,36.38) -- (557.53,42.91) ;
			\draw    (557.63,36.38) -- (564.23,42.91) ;
			\draw    (564.13,36.38) -- (570.73,42.91) ;
			\draw    (571.11,36.38) -- (577.7,42.91) ;
			\draw    (576.69,36.38) -- (583.29,42.91) ;
			\draw    (583.38,36.38) -- (589.98,42.91) ;
			\draw    (589.89,36.38) -- (596.49,42.91) ;
			\draw    (596.86,36.38) -- (603.46,42.91) ;
			\draw    (365.54,179.91) -- (372.14,186.43) ;
			\draw    (372.23,179.91) -- (378.83,186.43) ;
			\draw    (378.74,179.91) -- (385.33,186.43) ;
			\draw    (385.71,179.91) -- (392.31,186.43) ;
			\draw    (391.93,179.91) -- (398.53,186.43) ;
			\draw    (398.63,179.91) -- (405.22,186.43) ;
			\draw    (405.13,179.91) -- (411.73,186.43) ;
			\draw    (412.11,179.91) -- (418.7,186.43) ;
			\draw    (418.72,179.91) -- (425.32,186.43) ;
			\draw    (425.42,179.91) -- (432.01,186.43) ;
			\draw    (431.92,179.91) -- (438.52,186.43) ;
			\draw    (438.9,179.91) -- (445.49,186.43) ;
			\draw    (445.12,179.91) -- (451.72,186.43) ;
			\draw    (451.81,179.91) -- (458.41,186.43) ;
			\draw    (458.31,179.91) -- (464.91,186.43) ;
			\draw    (465.29,179.91) -- (471.89,186.43) ;
			\draw    (471.36,179.91) -- (477.95,186.43) ;
			\draw    (478.05,179.91) -- (484.65,186.43) ;
			\draw    (484.55,179.91) -- (491.15,186.43) ;
			\draw    (491.53,179.91) -- (498.13,186.43) ;
			\draw    (497.75,179.91) -- (504.35,186.43) ;
			\draw    (504.44,179.91) -- (511.04,186.43) ;
			\draw    (510.95,179.91) -- (517.55,186.43) ;
			\draw    (517.92,179.91) -- (524.52,186.43) ;
			\draw    (524.54,179.91) -- (531.14,186.43) ;
			\draw    (531.23,179.91) -- (537.83,186.43) ;
			\draw    (537.74,179.91) -- (544.34,186.43) ;
			\draw    (544.71,179.91) -- (551.31,186.43) ;
			\draw    (550.93,179.91) -- (557.53,186.43) ;
			\draw    (557.63,179.91) -- (564.23,186.43) ;
			\draw    (564.13,179.91) -- (570.73,186.43) ;
			\draw    (571.11,179.91) -- (577.7,186.43) ;
			\draw    (576.69,179.91) -- (583.29,186.43) ;
			\draw    (583.38,179.91) -- (589.98,186.43) ;
			\draw    (589.89,179.91) -- (596.49,186.43) ;
			\draw    (596.86,179.91) -- (603.46,186.43) ;
			\draw [color={rgb, 255:red, 65; green, 117; blue, 5 }  ,draw opacity=1 ][line width=1.5]    (365.89,126.33) .. controls (381.97,125.2) and (424.4,92.97) .. (461.07,81.63) .. controls (497.73,70.3) and (486.73,94.63) .. (451.4,116.63) .. controls (416.07,138.63) and (481.07,151.63) .. (506.73,149.3) .. controls (532.4,146.97) and (567.4,126.3) .. (602.46,126.54) ;

			\draw (49.11,2.4) node [anchor=north west][inner sep=0.75pt]    {$y$};
			\draw (260.84,132.32) node [anchor=north west][inner sep=0.75pt]    {$x\in \mathbb{T}^{d}$};
			\draw (8.22,32.87) node [anchor=north west][inner sep=0.75pt]    {$H^{-}$};
			\draw (-5,172.03) node [anchor=north west][inner sep=0.75pt]    {$-H^{+}$};
			\draw (144.07,93.11) node [anchor=north west][inner sep=0.75pt]  [color={rgb, 255:red, 208; green, 2; blue, 27 }  ,opacity=1 ]  {$\textcolor[rgb]{0.82,0.01,0.11}{\Sigma ( t) =\{y=\eta ( t,x)\}}$};
			\draw (57.44,61.94) node [anchor=north west][inner sep=0.75pt]  [color={rgb, 255:red, 74; green, 144; blue, 226 }  ,opacity=1 ]  {$\Omega ^{-}( t) \ $};
			\draw (56.71,143.63) node [anchor=north west][inner sep=0.75pt]  [color={rgb, 255:red, 245; green, 166; blue, 35 }  ,opacity=1 ]  {$\Omega ^{+}( t) \ $};
			\draw (374.11,2.4) node [anchor=north west][inner sep=0.75pt]    {$y$};
			\draw (585.84,133.65) node [anchor=north west][inner sep=0.75pt]    {$x\in \mathbb{T}$};
			\draw (333.22,32.87) node [anchor=north west][inner sep=0.75pt]    {$H^{-}$};
			\draw (320,172.03) node [anchor=north west][inner sep=0.75pt]    {$-H^{+}$};
			\draw (382.44,61.94) node [anchor=north west][inner sep=0.75pt]  [color={rgb, 255:red, 74; green, 144; blue, 226 }  ,opacity=1 ]  {$\Omega ^{-}( t) \ $};
			\draw (382.37,144.29) node [anchor=north west][inner sep=0.75pt]  [color={rgb, 255:red, 245; green, 166; blue, 35 }  ,opacity=1 ]  {$\Omega ^{+}( t) \ $};
			\draw (505.07,91.45) node [anchor=north west][inner sep=0.75pt]  [color={rgb, 255:red, 208; green, 2; blue, 27 }  ,opacity=1 ]  {$\textcolor[rgb]{0.82,0.01,0.11}{\Sigma }\textcolor[rgb]{0.82,0.01,0.11}{(}\textcolor[rgb]{0.82,0.01,0.11}{t}\textcolor[rgb]{0.82,0.01,0.11}{)}\textcolor[rgb]{0.82,0.01,0.11}{=}\textcolor[rgb]{0.82,0.01,0.11}{\{\gamma }\textcolor[rgb]{0.82,0.01,0.11}{(}\textcolor[rgb]{0.82,0.01,0.11}{t,s}\textcolor[rgb]{0.82,0.01,0.11}{)}\textcolor[rgb]{0.82,0.01,0.11}{\}}$};

		\end{tikzpicture}

		\caption{Graph case \& 2D overlapping case}
	\end{figure}
	
	Another one is the \textit{2D overlapping case}, where $d=1$ and the interface is a curve in 2D. In this case, we only consider the arc-length parameter $\gamma(t)\colon\R/L(t)\to\T\times]-H^+,H^-[$, where $L(t)$ denotes the total length of $\Sigma(t)$. To avoid self-intersection, one may add the \textit{chord-arc condition} (see (4) of Definition~\ref{def-intro-instab:RegSol} below).

	\subsection{Conservation laws and Hamiltonian structure}\label{subsect-intro:ConservationHamilton}
	
	For the two-phase water-wave system \eqref{eq-intro-ww:Main}, there are two frequently used conservation laws, the conservation of mass and energy. The mass is defined as
	\begin{equation}\label{eq-intro-ch:Mass}
		M := \int_\Sigma n\cdot(0,y) dS,
	\end{equation}
	where $dS$ denotes the surface element (see also \eqref{eq-pre-cons:MassAlt} for an alternative formulation). For simplicity, we shall assume that $M=0$ for all time (see (3) of Definition~\ref{def-intro-instab:RegSol} below). The total energy of the system \eqref{eq-intro-ww:Main} reads 
	\begin{equation}\label{eq-intro-ch:TotEnergy}
		E = E_k + E_p
	\end{equation}
	with kinetic energy $E_k$ and potential energy $E_p$ defined by
	\begin{align}
		E_k:=& \frac{1}{2} \sum_{\pm} \underline{\rho}^\pm \iint_{\Omega^\pm} |\nabla_{x,y}\phi^\pm|^2 dxdy, \label{eq-intro-ch:KinEnergy} \\
		E_p:=& \frac{Ag}{2}  \int_{\Sigma} n\cdot(0,y^2) dS + \frac{\sigma}{\rho^++\rho^-} \left[\mathrm{Area}(\Sigma(t)) - (2\pi)^d\right], \label{eq-intro-ch:PotEnergy}
	\end{align}
	where $\underline{\rho}^\pm = \rho^\pm/(\rho^++\rho^-)$ is the relative density, the constant
	\begin{equation}\label{eq-intro-ch:AtwoodNum}
		A := \frac{\rho^+-\rho^-}{\rho^++\rho^-}
	\end{equation}
	is \textit{Atwood number}, and $\mathrm{Area}(\Sigma(t))$ is the area of the hypersurface $\Sigma(t)$. This energy is preserved in time, while it is not positively definite when the Atwood number $A<0$. A detailed justification of these conservation laws will be left to Section~\ref{subsect-pre:Conservation}.
	
	Consider the unified trace of scalar potential,
	\begin{equation}\label{eq-intro-ch:DefPsi}
		\psi := \underline{\rho}^+ \phi^+|_{\Sigma} - \underline{\rho}^- \phi^-|_{\Sigma}.
	\end{equation}
	If, in addition, the interface $\Sigma(t)$ can be represented as $y=\eta(t,x)$ for some function $\eta$, the energy $E$ will be a nonlinear functional of $(\eta,\psi)$ and the system \eqref{eq-intro-ww:Main} admits the following Hamiltonian formulation,
	\begin{equation}\label{eq-intro-ch:Hamilton}
		\eta_t = \frac{\delta E}{\delta \psi},\ \ \psi_t = -\frac{\delta E}{\delta \eta}.
	\end{equation}
	A formal derivation of this can be found in \cite{kuznetsov1995nonlinear,benjamin1997reappraisal} for graph case and in \cite{craig2017hamiltonian} for overlapping case. In view of this (formal) Hamiltonian structure, one may expect to establish a \textit{virial theorem} for the system \eqref{eq-intro-ww:Main}, which will be discussed in the next subsection.

	\subsection{Virial theorem}\label{subsect-intro:Virial}
	
	The \textit{virial theorem} describes the phenomenon that, in many Hamiltonian systems, the kinetic and potential energy are equal in the time-averaged sense. Enunciated by Clausius \cite{clausius1870xvi} in 1870, it quickly became a powerful tool in physics \cite{zwicky1933rotverschiebung,rayleigh1905xlii,einstein1922grundlage}. The connection between virial theorem and numerous important models in classical physics has also been discovered in recent decades \cite{hanson1995virial,georgescu1999virial,shore2012introduction,alazard2023virial}. In addition, the idea of this theorem could also be applied to some topics in biology \cite{podio2019virial} and economics \cite{liorsdottir2023virial}. 
	
	In mathematics, virial theorem serves as an important method in the blow-up problem of PDEs. In the 1970s, Zakharov-Sobolev-Synakh \cite{zakharov1971character} and Glassey \cite{glassey1977blowing} proved the virial theorem for nonlinear Schr\"odinger equation and deduced the existence of singular solution in defocusing settings. Such technique is also used by Merle \cite{merle1996blow}, Sideris \cite{sideris1985formation}, Levine \cite{levine1974instability,levine1990role}, Keel-Tao \cite{keel1999small}, and Kenig-Merle \cite{kenig2008global}. Recently, Alazard-Zuily \cite{alazard2023virial} proved the first virial theorem for water-wave equation and established the Rayleigh-Taylor instability for a large class of initial data, which will be discussed in the next part. 
	
	Additionally, the virial theorem is often expressed as \textit{equipartition of energy} which asserts that the kinetic and potential energy are equal, asymptotically in time. This phenomenon itself has also attracted lots of attentions in the mathematics and physics community. We refer to \cite{costa1981energy} for the case of wave equation and Maxwell's equation, \cite{glassey1973asymptotic} for nonlinear wave equation, \cite{dassios1984dissipation} for thermoelasticity, and \cite{vega2008equipartition} for more results on nonlinear wave equation.
	
	Inspired by \cite{alazard2023virial}, we shall prove in this article a virial theorem for the two-phase water-wave problem \eqref{eq-intro-ww:Main}. 
	\begin{theorem}\label{thm-intro-viri:Virial}
		Let $(\gamma, \phi^\pm,P^\pm)$ be a regular solution to \eqref{eq-intro-ww:Main} on $[0,T[$ (see Definition~\ref{def-intro-instab:RegSol} below) for some $T\in]0,+\infty]$. Then the following identity holds for all time $t\in [0,T[$,
		\begin{equation}\label{eq-intro-viri:Main}
			\frac{1}{2}\frac{d}{dt} \int_{\Sigma}n\cdot(0,y\psi) dS = \tilde{E}_k - E_p + R,
		\end{equation}
		where $\psi$ is defined in \eqref{eq-intro-ch:DefPsi}, $\tilde{E}_k$ is the \textit{modified kinetic energy}
		\begin{equation}\label{eq-intro-viri:ModKinEnergy}
			\tilde{E}_k:= \sum_{\pm} \underline{\rho}^\pm \iint_{\Omega^\pm} \left( \frac{1}{4}|\nabla_{x}\phi^\pm|^2 + \frac{3}{4}|\partial_y\phi^\pm|^2 \right)  dxdy,
		\end{equation}
		$E_p$ is the potential energy defined in \eqref{eq-intro-ch:PotEnergy}, and $R$ is a non-negative remainder.
		
		The remainder $R$ reads
		\begin{equation}\label{eq-intro-viri:Rem}
			R = R_{b}^+ + R_{b}^- + R_s,
		\end{equation}
		with
		\begin{align}
			R_{b}^\pm =& \left\{\begin{array}{ll}
				\frac{1}{4} \underline{\rho}^\pm H^\pm \int_{\Gamma^\pm}|\nabla_x\phi^\pm|^2dx, & \text{when }H^\pm<+\infty, \\[1ex]
				0, & \text{when }H^\pm=+\infty,
			\end{array}\right. \label{eq-intro-viri:RemBottom}\\
			R_s =& \frac{\sigma}{2(\rho^+ + \rho^-)} \left[ \mathrm{Area}(\Sigma) + \int_{\Sigma}|n\cdot(0,1)|^2 dS - 2(2\pi)^d \right].\label{eq-intro-viri:RemSurfaceTension}
		\end{align}
		Here $R_b^\pm$ and $R_s$ are both non-negative.
	\end{theorem}
	
	Note that, in \cite{alazard2023virial} where the authors focus on the one-phase graph case (namely, $\rho^+=1$ and $\rho^-=0$), the integral on the left-hand side of \eqref{eq-intro-viri:Main} is the inner product of $\eta$ and $\psi$, which are variables in the Hamiltonian formulation \eqref{eq-intro-ch:Hamilton}. Our choice here is a general version of the same quantity.
	
	Under the condition $A\le 0$ together with $E\neq 0$ when $\sigma=0$ and $E<0$ when $\sigma>0$ (see Corollary~\ref{cor-ins-viri:Main}), the right hand side of \eqref{eq-intro-viri:Main} can be bounded from below by a positive time-independent constant. Then this theorem implies that the integral on the left-hand of \eqref{eq-intro-viri:Main} grows at least faster than $O(t)$, so as its upper bounds. Thus, we will attempt to control this integral by a function of some geometrical quantities (for instance, slope and curvature) and deduce their growth. The detailed results will be presented in the next paragraph.

	\subsection{Instability of the system}\label{subsect-intro:Instab}
	
	\subsubsection{Well-posedness of the problem}
	
	In the two-phase water-wave problem \eqref{eq-intro-ww:Main}, if we assume that the \textit{Atwood number} $A$ defined by \eqref{eq-intro-ch:AtwoodNum} is negative, meaning that the denser fluid is placed above the other one, the system will become highly unstable, which is known as the \textit{Rayleigh-Taylor Instability} (RTI), due to the pioneer work of Rayleigh \cite{rayleigh1883investigation} and Taylor \cite{taylor1950instability}. A comprehensive physical description of this problem can be found in \cite{kull1991theory}. In the degenerate case $A=0$, or equivalently $\rho^+=\rho^-$, another instability arises from the jump of tangential velocity instead of the gravity, known as the \textit{Kelvin-Helmholtz Instability} (KHI), which was firstly studied by Kelvin \cite{thomson1871xlvi} and Helmholtz \cite{helmholtz1868xliii}. Let us first review the well-posedness of RTI and KHI problems.
	
	In the absence of surface tension, the local well-posedness theory can be established in the space of analytical functions. For example, in \cite{sulem1981finite,sulem1985finite}, the authors used the Cauchy-Kowalewski theorem to prove that if the initial data could be extended to a complex strip of width $r_0>0$, then there exists a unique solution in the space formed by analytic functions in a strip of width $r\in ]0,r_0[$, and the existence time is of size $r_0-r$. This result has only been proved for small initial data, where the interface has small gradient and is flat at infinity. This smallness condition is not essential in 2D torus case due to Xie \cite{xie2007existence} by pulling the system to unit disc via conformal mapping.
	
	Indeed, the analytic framework is optimal, since the solution to KHI or RTI problem is either excessively rough or analytic. More precisely, for 2D KHI problem, Lebeau proved in \cite{lebeau2002regularite} that, once the interface has H\"older regularity $C^{1,\epsilon}$ with $\epsilon>0$, the solution must be analytic in time and space. This result was also proved independently by Wu in \cite{wu2006mathematical}, where the interface is assumed to have only Lipschitz regularity. And for 2D RTI problem, due to \cite{kamotski2005rayleigh}, H\"older regularity $C^{1,\epsilon}$ of interface guarantees that the solution is smooth, or analytic in space with the additional assumption that the jump at interface is non-vanishing. Remark that the minimal regularity assumed for the interface should be Lipschitz in order to give a proper sense to the system \eqref{eq-intro-ww:Main}. In addition, if the irrotational condition is eliminated, Ebin proved in \cite{ebin1988ill} that both RTI and KHI problems are ill-posed (Hadamard sense) in Sobolev spaces. A similar result was also proved by Guo-Tice \cite{guo2011compressible} for compressible fluids. 
	
	In the presence of surface tension, most of the mathematical researches are dedicated to KHI, or equivalently the \textit{vortex sheet problem}. In this case, the surface tension can be shown to stabilize the system and better local well-posedness results could be expected, compared with the case without surface tension. For example, in \cite{iguchi1997twophase} the authors proved the local well-posedness in high order Sobolev spaces for small initial data where the interface is almost flat and the velocity is close to zero. The smallness condition is then eliminated by Ambrose \cite{ambrose2003wp} in 2D and by Ambrose-Masmoudi \cite{ambrose2007wp} in 3D. We also refer to \cite{shatah2008apriori,cheng2008motion} for similar results.
	
	For RTI with surface tension, seldom results are obtained for inviscid fluids. However, when the viscosity is considered, the well-posedness theory will be available, as in the case of water-wave problem \cite{guo2013local, guo2013decay, guo2013almost}. For incompressible viscous fluids, in \cite{wang2014viscous}, the authors proved local well-posedness in high order Sobolev spaces for small data and their almost exponential decay. In compressible viscous setting, it is proved in \cite{jang2016compressible1,jang2016compressible2,jang2016compressible3} that the problem is locally well-posed for small initial data in high order Sobolev spaces and is exponentially stable near equilibrium when the surface tension is large enough.
	
	In this article, we will not concern the well-posedness problem, but concentrate on the characterization of instability for the system \eqref{eq-intro-ww:Main}. As discussed above, it is reasonable to assume that there exist non-trivial solutions to the main system \eqref{eq-intro-ww:Main} and that these solutions are smooth enough (see (1) of Definition~\ref{def-intro-instab:RegSol} below), for the sake of simplicity. Before stating our main results, we precise the solutions to be considered throughout this article.

	\begin{definition}[Regular solutions] \label{def-intro-instab:RegSol}
		Given $T\in]0,+\infty]$, $g,\sigma\in[0,+\infty[$, $H^\pm\in]0,+\infty]$, and $\rho^\pm\in[0,+\infty[$ with $\rho^++\rho^->0$, we say $(\gamma,\phi^\pm,P^\pm)$ is a regular solution to the system \eqref{eq-intro-ww:Main} on $[0,T[$, if it solves \eqref{eq-intro-ww:Main} on the time interval $[0,T[$ with the following conditions:
		\begin{enumerate}[(1)]
			\item (Smoothness) The parametrization $\gamma$, the scalar potential $\phi^\pm$, and the pressure $P^\pm$,
			$$
			\gamma\colon[0,T[\times\T^d\to\Omega, \quad \phi^\pm,P^\pm\colon\{(t,x,y):(x,y)\in\Omega(t)\}\to\R,
			$$
			are all $C^1$ in time. In space variable, $\gamma(t)$ is $C^2$, $\phi^\pm(t)$ is $C^2$ up to boundary, and $P^\pm(t)$ is $C^1$ up to boundary. 
			\item (Finite energy) The kinetic energy $E_k = E_k(t)$ defined by \eqref{eq-intro-ch:KinEnergy} is finite for all time $t\in[0,T[$ with $E_k \in L^\infty_\mathrm{loc}([0,T[)$.
			\item (Zero mass) The mass $M$ defined by \eqref{eq-intro-ch:Mass} is zero for all time $t\in[0,T[$.
			\item (No self-intersection) For all time $t\in[0,T[$, the parametrization $\gamma(t)$ is a diffeomorphism from $\T^d$ to $\Sigma(t)$. In particular, for 2D case (namely $d=1$), we take $\gamma(t)$ as the arc-length parametrization and add the chord-arc condition: there exists $0<c_0\le 1$, such that, for all $s_1,s_2\in\R/L(t)$ and all time $t\in[0,T[$,
			\begin{equation}\label{eq-intro-instab:ChordArc}
				c_0 |s_1-s_2| \le |\gamma(t,s_1)-\gamma(t,s_2)| \le |s_1-s_2|,
			\end{equation}
			where $L(t)$ is the length of $\Sigma(t)$ and the absolute value $|s_1-s_2|$ should be understood as the distance in $\R/L(t)$.
			\item (Distance between free and solid boundaries) If $H^+<+\infty$ (resp. $H^-<+\infty$), there exists $d_0>0$ such that, for all time $t\in[0,T[$, the distance between the free boundary $\Sigma(t)$ and the solid boundary $\Gamma^+$ (resp. $\Gamma^-$) is larger than $d_0$. Moreover, when $H^+ = H^- = +\infty$, we take $d_0=+\infty$.
			\item (Growth of pressure) If $H^+=+\infty$ (resp. $H^-=+\infty$), the pressure $P^+$ (resp. $P^-$) has at most polynomial growth for all time $t\in[0,T[$. Namely, for arbitrary $0<T'<T$, there exists $N\in\N$ such that
			\begin{equation}\label{eq-intro-instab:GrowPress}
				\sup_{x\in\T^d,\ t\in[0,T']} \left| P^\pm(t,x,y) \right| \lesssim (1+|y|)^N.
			\end{equation}
		\end{enumerate}
		Moreover, $(\gamma, \phi^\pm,P^\pm)$ is said to be global-in-time if $T=+\infty$.
	\end{definition}
	
	\begin{remark}\label{rmk-intro-instab:DecayScalPot}
		In infinite depth case where $H^+=+\infty$ (resp. $H^-=+\infty$), one could deduce from the conditions (1) and (2), that the derivatives of the scalar potential $\phi^+$ (resp. $\phi^-$) is $C^\infty$ in space and decays rapidly at infinity. Rigorously, for arbitrary sub-interval $[0,T']\subset[0,T[$, we have, for all $k\in\N$ and $\alpha\in\N^{d+1}$ with $|\alpha|\ge 1$,
		\begin{equation}\label{eq-intro-instab:DecayScalPot}
			\begin{aligned}
				&\lim_{y\to -\infty} \sup_{x\in\T^d,\ t\in[0,T']} \left| y^k \nabla_{x,y}^\alpha\phi^+(t,x,y) \right| = 0,\\
				\text{resp. }&\lim_{y\to +\infty} \sup_{x\in\T^d,\ t\in[0,T']} \left| y^k \nabla_{x,y}^\alpha\phi^-(t,x,y) \right| = 0.
			\end{aligned}
		\end{equation}
		The proof of this decay is based on classical elliptic theories and will be detailed in Section~\ref{subsect-pre:Decay}. 
	\end{remark}
	
	\begin{remark}\label{rmk-intro-instab:Splash}
		The conditions (4) in Definition~\ref{def-intro-instab:RegSol} is assumed so that the system \eqref{eq-intro-ww:Main} is valid. However, it is unclear whether conditions (4) can hold for large time. In fact, the self-intersection, also known as \textit{splash singularity} has already been proved for 2D water-wave problem ($\rho^+=1$ and $\rho^-=0$), cf. \cite{castro2012splash}, while the corresponding result for two-phase water-wave is still open. 
	\end{remark}
	
	\begin{remark}\label{rmk-intro-instab:Mixing}
		The conditions (5) in Definition~\ref{def-intro-instab:RegSol} is also proposed for the validation of the system \eqref{eq-intro-ww:Main}. One should note that it fails in finite time for RTI with finite depths. Indeed, it is well-known that the free boundary $\Sigma(t)$ is contained in a strip $\T^d\times[-\alpha^+(t),\alpha^-(t)]$, known as mixing area, with upper and lower heights $\alpha_\pm(t)$ of $\Sigma(t)$ growing like $O(t^2)$, cf. \cite{chandrasekhar2013hydrodynamic,kalinin2024scale,gebhard2021newapproach}. Therefore, the growth of mixing area could invalidate the system \eqref{eq-intro-ww:Main}, which is not the subject of this article.
	\end{remark}
	
	\begin{remark}\label{rmk-intro-instab:Press}
		The purpose of the condition (6) in Definition~\ref{def-intro-instab:RegSol} is to avoid exponential growth of $P^\pm$ at infinite. One can see from the first equation of \eqref{eq-intro-ww:Main} that $P^\pm$ solves the elliptic equation
		\begin{equation}\label{eq-intro-instab:EqPress}
			\Delta_{x,y}P^\pm = - \frac{\rho^\pm}{2}\Delta_{x,y}\left(|\nabla_{x,y}\phi^\pm|^2\right)
		\end{equation}
		at least in the sense of distribution, where the right-hand side decays rapidly at infinity (see \eqref{eq-intro-instab:DecayScalPot}). Through separation of variables, $P^\pm$ can be written as the sum of a term with at most polynomial growth in $y$ and an exponentially increasing (in $y$) harmonic function. With the condition (6) in Definition~\ref{def-intro-instab:RegSol}, this extra harmonic function will not appear.
	\end{remark}
	
	\subsubsection{Stability criteria and the description of instability}
	
	In the one-phase case ($\rho^+>0$ and $\rho^-=0$), a well-known stability criteria is \textit{Taylor sign condition} : 
	$$
	-\partial_y P^+|_{\Sigma} >0,
	$$
	which is rigorously justified by Wu in 1990s \cite{wu1997well, wu1999well}. For two-phase fluid, an analogue is given by Kelvin for the linearized system : 
	$$
	g(\rho^+ - \rho^-) > \frac{1}{4\sigma}\left(\frac{\rho^+\rho^-}{\rho^++\rho^-}\right)^2 |\BBra{u}|^4,
	$$
	where $\BBra{u}$ denotes the jump of the horizontal velocity at the interface (see for example Ch.X and Ch.XI of \cite{chandrasekhar2013hydrodynamic}). This stability condition is then generalized to the nonlinear problem (graph case) by Lannes \cite{lannes2013stability} :
	$$
	\BBra{-\partial_y P} > \frac{1}{4\sigma}\left(\frac{\rho^+\rho^-}{\rho^++\rho^-}\right)^2 \mathfrak{c}(\eta) \|\BBra{V}\|_{L^\infty}^4,
	$$
	where $\BBra{-\partial_y P} := -\partial_y P^+|_{\Sigma} + \partial_y P^-|_{\Sigma} $ is the jump of $y$-derivative of pressure at the interface, $\mathfrak{c}(\eta)>0$ is a constant depending on the shape of the interface, $\BBra{V}$ is the jump of the horizontal velocity at the interface.
	
	According to these stability criteria, two potential sources of instability can be identified: gravity and the velocity jump at the interface. Evidently, the case with a negative Atwood number $A < 0$—or equivalently $\rho^+ < \rho^-$—corresponds to an unstable regime, known as the Rayleigh–Taylor instability. In contrast, the Kelvin–Helmholtz instability is driven by the discontinuity in velocity across the interface. When there is no density jump, i.e., $\rho^+ = \rho^-$, the RTI is absent, and the only remaining mechanism of instability is the velocity jump (KHI). It is also noteworthy that even when the Atwood number is positive ($A > 0$, i.e. $\rho^+ > \rho^-$), a sufficiently large velocity jump can still lead to instability, as indicated by the aforementioned stability criteria.
	
	In this article, we focus on the instable regime $A<0$ (namely, $\rho^+<\rho^-$) and attempt to obtain a quantitative description of the instability. From the linear analysis (see for example the summary \cite{bardos2010mathematics}), one may expect an exponential growth in short time. In fact, the existing results, mainly based on growing mode or maximal regularity \cite{hwang2003dynamic,guo2010linear,wang2012viscous,jang2016compressible2,jiang2023rayleigh,pruess2010rayleigh,wilke2017rayleigh} manage to find some data with such growth. Even if these methods allow us to trait more general contexts (compressible, rotational, viscous, etc), they do lay emphasis on the comparison between linear and nonlinear systems, from which one can only deduce the existence of unstable perturbation near zero in short time. Due to the physical nature of this problem, one should expect the blow-up (or, at least, rapid growth) for most regular solutions, which could be extremely challenging due to the full nonlinearity of the system. As indicated in \cite{alazard2023virial}, a possible solution could be the virial identity (see Theorem~\ref{thm-intro-viri:Virial}), which implies the following main results.
	
	\subsubsection{Main results}
	
	For graph case, we shall prove that the slope of the interface $\Sigma(t)$ grows faster than $O(t^{2/5})$ (or $O(t^{2/3})$ when $\rho^+=\rho^-$). Rigorously, we have
	\begin{theorem}[Instability of graph case]\label{thm-intro-instab:Graph}
		Consider $d\ge 1$, $H^\pm\in]0,+\infty]$, $\rho^+\le\rho^-$, $\sigma\ge0$, and $T\in]0,+\infty]$. Let $(\gamma, \phi^\pm,P^\pm)$ be any regular solution to \eqref{eq-intro-ww:Main} on $[0,T[$ (see Definition~\ref{def-intro-instab:RegSol}) with $\gamma(t,x)=(x,\eta(t,x))$ for some real-valued function $\eta$. We assume that the total energy $E$ defined by \eqref{eq-intro-ch:TotEnergy} satisfies 
		\begin{equation}\label{eq-intro-instab:CondTotEnergy}
			E\neq 0\text{ when }\sigma=0,\text{ and }E<0\text{ when }\sigma>0.
		\end{equation}
		Then there exists a constant $C>0$ depending on $d_0$ (see (5) in Definition~\ref{def-intro-instab:RegSol}), such that, for all time $t\in[0,T[$,
		\begin{equation}\label{eq-intro-instab:Graph}
			|E| t + \left.\int_{\T^d} \eta\psi dx\right|_{t=0}  \leqslant C \|\nabla_x\eta\|_{L^\infty} \sqrt{1+ \|\nabla_x\eta\|_{L^\infty}} \sqrt{|E| + |A|g\|\nabla_x\eta\|_{L^\infty}^2},
		\end{equation}
		where $\psi$ is defined in \eqref{eq-intro-ch:DefPsi} and $A\le 0$ is the Atwood number (see \eqref{eq-intro-ch:AtwoodNum}).
	\end{theorem}
	
	\begin{remark}
		When $\rho^+ < \rho^-$, the algebric growth observed in~\eqref{eq-intro-instab:Graph} may be due to Rayleigh–Taylor instability, as the condition $\rho^+ < \rho^-$ is critical in the proof. However, this explanation does not apply in the case $\rho^+ = \rho^-$. Recall that in the formulation of the problem, zero velocity jump across the interface is permitted, which precludes the direct connection with Kelvin–Helmholtz instability. In fact, as shown in the formulation via vorticity in~\cite{bardos2010mathematics}, the system is hyperbolic with degeneracy when $\rho^+ = \rho^-$, which may naturally lead to algebric growth in time.
	\end{remark}
	
	Among the conditions of Theorem~\ref{thm-intro-instab:Graph}, except for the reasonable assumptions listed in Definition~\ref{def-intro-instab:RegSol}, the only non-trivial one is \eqref{eq-intro-instab:CondTotEnergy} about the total energy, depending only on the initial data. Without surface tension, $E\neq0$ contains already a large class of solutions. However, in the case with surface tension ($\sigma>0$), we have to make the assumption $E<0$ so that the right-hand side of the virial identity \eqref{eq-intro-viri:Main} is positive. Note that $E<0$ can occur only when $A<0$ (see Remark~\ref{rmk-pre-cons:SignPotEnergy}) and it can be verified only by a few solutions such as those with initial data $\phi^\pm(0,x,y)=0$, $\eta(0,x)=M\sin(x)$, where $M\in\R$ is large enough. The characterization of instability for $E>0$, $\sigma>0$ is still open.
	
	For 2D overlapping case, instead of the slope, we shall consider the curvature of the interface $\Sigma(t)$. In a similar manner, the $L^\infty$-norm of the curvature can be shown to grow faster than $O(t^{2/9})$ (or $O(t^{1/3})$ when $\rho^+=\rho^-$), provided that the depths are finite.
	\begin{theorem}[Instability of 2D overlapping case]\label{thm-intro-instab:NonGraph}
		Consider $d=1$, $H^\pm<+\infty$, $\rho^+\le\rho^-$, $\sigma\ge0$, and $T\in]0,+\infty]$. Let $(\gamma, \phi^\pm,P^\pm)$ be any regular solution to \eqref{eq-intro-ww:Main} on $[0,T[$ (see Definition~\ref{def-intro-instab:RegSol}) with $\gamma$ being the arc-length parametrization. 
		
		If the total energy $E$ defined by \eqref{eq-intro-ch:TotEnergy} satisfies the condition \eqref{eq-intro-instab:CondTotEnergy}, then there exists a universal constant $N_0>0$ and another constant $C>0$ depending on $H^\pm$, such that, for all time $t\in[0,T[$,
		\begin{equation}\label{eq-intro-instab:NonGraph}
			|E| t + \left.\int_{\Sigma} n\cdot(0,y\psi) dS\right|_{t=0}  \leqslant C \epsilon(t)^{-3} \sqrt{|E|+|A|g\epsilon(t)^{-3}},
		\end{equation}
		with $\epsilon(t)>0$ defined by
		\begin{equation}\label{eq-intro-instab:DefEpsilon}
			\epsilon(t):= \min\left( \frac{c_0}{N_0(\|\kappa(t)\|_{L^\infty}+1)}, d_0 \right),
		\end{equation}
		where $\psi$ is defined in \eqref{eq-intro-ch:DefPsi}, $\kappa$ is the curvature, $c_0$ is the constant from the chord-arc condition \eqref{eq-intro-instab:ChordArc}, and $d_0$ is the constant appearing in the condition (5) of Definition~\ref{def-intro-instab:RegSol}.
	\end{theorem}
	
	\begin{remark}
		In infinite depth case $H^\pm=+\infty$, if there exists time-independent constants $0<h^\pm<+\infty$ such that the interface $\Sigma(t)$ lies in the strip $\{-h^+<y<h^-\}$ for all time $t\in [0,T]$, we can also prove the inequality \eqref{eq-intro-instab:NonGraph} by using the same argument with $H^\pm$ replaced by $h^\pm$. Note that, as mentioned in Remark~\ref{rmk-intro-instab:Mixing}, the upper and lower heights of $\Sigma(t)$ have $O(t^2)$ growth. Thus, the assumption $\Sigma(t)\subset \{-h^+<y<h^-\}$ must fail in finite time.
	\end{remark}
	
	Note that if the maximal time of existence $T$ is large enough, the right-hand side of \eqref{eq-intro-instab:NonGraph} will be large enough so that $\epsilon^{-1} > d_0^{-1}$, which implies the growth of $\|\kappa(t)\|_{L^\infty}$. This result should be understood as: if the self-intersection and the intersection between free and solid boundaries never occur, then the curvature $\|\kappa(t)\|_{L^\infty}$ grows faster than $O(t^{2/9})$ (or $O(t^{1/3})$ for KHI). That is to say, in general, there are at least three possible singularities in large time, 
	\begin{enumerate}[(1)]
		\item The curvature $\|\kappa\|_{L^\infty}$ tends to infinity (faster than $O(t^{2/9})$ for RTI and $O(t^{1/3})$ for KHI);
		\item There exists self-intersection (see Remark~\ref{rmk-intro-instab:Splash});
		\item The free boundary $\Sigma(t)$ intersects with the solid boundaries $\Gamma^\pm$ (see Remark~\ref{rmk-intro-instab:Mixing}).
	\end{enumerate}
	
	The proof of Theorem~\ref{thm-intro-instab:Graph} and~\ref{thm-intro-instab:NonGraph} are based on the virial identity \eqref{eq-intro-viri:Main}. Noticing that the right-hand side of \eqref{eq-intro-viri:Main} can be bounded from below by $|E|/2$ (see Corollary~\ref{cor-ins-viri:Main}), it suffices to control the integral of $n\cdot(0,y\psi)$ on $\Sigma(t)$, where the main difficulty is to bound $\psi$ by quantities depending only on slope or curvature. Since $\psi$ (see \eqref{eq-intro-ch:DefPsi}) is defined using the trace of scalar potential $\phi^\pm$, one may expect to apply the trace estimate. Namely, the $L^2$-norm of $\psi$ can be bounded by the $L^2$-norm of $\nabla_{x,y}\phi^\pm$, or equivalently, the kinematic energy $E_k$ (see \eqref{eq-intro-ch:KinEnergy}). Then the conservation of energy allows us to write $E_k = E - E_p$ with $E_p$ depending only on the geometry of the interface $\Sigma(t)$, which can be controlled by slope or curvature via some direct calculus.
	
	\subsection{Plan of the article and conventions}\label{subsect-intro:Notations}
	
	In Section~\ref{Sect:Pre}, we shall clarify the decay of scalar potential (in infinite depth case) and a proof of conservation laws will be given. Then the virial theorem~\ref{thm-intro-viri:Virial} will be proved in Section~\ref{Sect:viri}. Based on this, we conclude the instability results, Theorem~\ref{thm-intro-instab:Graph} and~\ref{thm-intro-instab:NonGraph}, in Section~\ref{Sect:instab}, where the graph case and 2D overlapping case are treated separately. In Appendix~\ref{App:nor-geo}, we will show that the normal geodesic coordinate realizes a global diffeomorphism, which is useful in the proof of Theorem~\ref{thm-intro-instab:NonGraph}.
	
	Throughout this article, we shall frequently use the following notations and conventions:
	\begin{enumerate}
		\item $\T$ represents the periodic domain $\R/(2\pi)$.
		\item For non-negative quantities $A,B$ and some parameter $\alpha$, we write $A\lesssim_\alpha B$, if there exists a constant $C_\alpha>0$ depending on $\alpha$ such that $A \le C_\alpha B$. Furthermore, we write $A\sim_\alpha B$ when $A\lesssim_\alpha B$ and $B\lesssim_\alpha A$.
		\item For non-negative quantities $A,B$, we write $A\ll B$ if there exists a large universal constant $C>0$ such that $A \le CB$.
		\item In infinite depth case, the integral of a function $f$ on $\Omega^\pm$ should be understood as the limit
		\begin{equation*}
			\iint_{\Omega^\pm}fdxdy = \lim_{h^\pm\to+\infty} \iint_{\Omega^\pm\cap\{\mp y<h^\pm\}}fdxdy.
		\end{equation*}
		The functions $f$ to be studied in this article decay to zero as $y\to\pm\infty$, uniformly in $x\in\T^d$. The precised behavior of $f$ at infinity will be clarified in the context.
	\end{enumerate}

	\medskip
	
	\noindent\textbf{Acknowledgements.} The author would like to warmly thank the referees for numerous suggestions which help to greatly improve the presentation of these results.

	\section{Preliminaries}\label{Sect:Pre}
	
	The aim of this section is to rigorously prove some results well-known in physics, including the decay of scalar potential $\phi^\pm$ and the conservation of mass and energy.

	\subsection{Decay of scalar potential}\label{subsect-pre:Decay}
	
	In this part, we clarify how to deduce the rapid decay \eqref{eq-intro-instab:DecayScalPot} of the scalar potentials in infinite depth cases. For simplicity, we will focus on the lower fluid domain $\Omega^+$ with $H^+=+\infty$ (the argument for the upper domain is the same). To conclude the decay \eqref{eq-intro-instab:DecayScalPot}, we choose a large depth $-h^+$ that the free boundary $\Sigma$ does not reach.  By using the data at $y=-h^+$ and the condition $E_k<+\infty$, we will find an explicit formula for $\phi^+$, which simply gives the desired decay. The detailed result is stated in the following lemma.
	
	\begin{lemma}\label{lem-pre-decay:Main}
		Let $(\gamma,\phi^\pm,P^\pm)$ be a regular solution to \eqref{eq-intro-ww:Main} on the time interval $[0,T[$. If the depth $H^+=+\infty$ (resp. $H^-=+\infty$), then, for all time $t\in[0,T[$, $\phi^\pm(t)$ is $C^\infty$ inside $\Omega^\pm(t)$. Moreover, for all $0<T'<T$, $k\in\N$, and $\alpha\in\N^{d+1}$ with $|\alpha|\ge 1$, there holds
		\begin{equation}\label{eq-pre-decay:Main}
			\begin{aligned}
				&\lim_{y\to-\infty} \sup_{x\in\T^d,\ t\in[0,T']} \left| y^k \nabla_{x,y}^\alpha\phi^+(t,x,y) \right| = 0, \\
				\text{resp. }&\lim_{y\to+\infty} \sup_{x\in\T^d,\ t\in[0,T']} \left| y^k \nabla_{x,y}^\alpha\phi^-(t,x,y) \right| = 0.
			\end{aligned}
		\end{equation}
	\end{lemma}
	\begin{proof}
		The smoothness of $\phi^\pm(t)$ is a consequence of elliptic regularity since $\Delta_{x,y}\phi^\pm=0$. For simplicity, in what follows, we focus on the decay of $\phi^+$, while $\phi^-$ can be treated in the same way. 
		
		Since the parametrization $\gamma$ is continuous on $[0,T']\times\T^d$, there exists $h^+>0$ such that $\T^d\times]-\infty,-h^+]$ is contained in $\Omega^\pm(t)$, for all $t\in[0,T']$. Since $\phi^+$ is smooth in space variables, we defined two functions $f,g \in C^\infty(\T^d)$ by $f(x)=\phi^+(x,-h^+)$ and $g(x) = \partial_y\phi^+(x,-h^+)$. Here the time dependence is omitted for simplicity, since all the arguments below are uniform in time $t\in[0,T']$. The smoothness of $f,g$ guarantees that their Fourier coefficients have arbitrary polynomial decay. Namely, for all $N\in\N$,
		\begin{equation}\label{eq-pre-decay:DecayFourierCoeef}
			|\hat{f}(\xi)| + |\hat{g}(\xi)| \lesssim_N (1+|\xi|)^{-N}.
		\end{equation}
		Thanks to $\Delta_{x,y}\phi^+=0$ (the second equation in \eqref{eq-intro-ww:Main}), the standard Fourier method gives the following explicit formula.
		\begin{equation*}
			\begin{aligned}
				\hat{\phi}^+(\xi,y) =& \frac{1}{2}e^{(y+h^+)|\xi|}\left( \hat{f}(\xi) + |\xi|^{-1}\hat{g}(\xi) \right) + \frac{1}{2}e^{-(y+h^+)|\xi|}\left( \hat{f}(\xi) - |\xi|^{-1}\hat{g}(\xi) \right), \quad\text{if }\xi\neq 0 \\
				\hat{\phi}^+(0,y) =& \hat{g}(0)(y+h^+) + \hat{f}(0),
			\end{aligned}
		\end{equation*}
		where $\hat{\phi}^+(\xi,y)$ is the Fourier coefficient of $\phi^+(y)$ at frequency $\xi$. We claim that $\hat{g}(\xi)=|\xi|\hat{f}(\xi)$ for all $\xi\in\Z^d$. Otherwise, the $L^2(\T^d)$-norm of $\nabla_{x,y}\phi^+(\cdot,y)$ admits a strictly positive lower bound as $y\to-\infty$, contradicting the condition (2) in Definition~\ref{def-intro-instab:RegSol} which implies $\nabla_{x,y}\phi^+\in L^2(\Omega^+)$. Consequently, we have
		\begin{equation*}
			\hat{\phi}^+(\xi,y) = \frac{1}{2}e^{(y+h^+)|\xi|}\left( \hat{f}(\xi) + |\xi|^{-1}\hat{g}(\xi) \right), \quad\text{if }\xi\neq 0,\quad \text{and } \hat{\phi}^+(0,y) = \hat{f}(0).
		\end{equation*}
		Based on this, we compute, for all $k\in\N$ and $\alpha\in\N^{d+1}$ with $|\alpha|\ge 1$,
		\begin{equation*}
			\left| y^k \nabla_{x,y}^\alpha\phi^+(x,y) \right| \lesssim |y|^k \sum_{\xi\in\Z^d} e^{(y+h^+)|\xi|}|\xi|^{|\alpha|-1} \left||\xi|\hat{f}(\xi) + \hat{g}(\xi)\right|.
		\end{equation*}
		Here the choice $|\alpha|\ge 1$ ensures that there is no contribution from $\xi=0$ in the sum on the right-hand side. Then due to the rapid decay \eqref{eq-pre-decay:DecayFourierCoeef} of $\hat{f}(\xi),\hat{g}(\xi)$, the right-hand side, independent of $x$, converges to zero as $y\to-\infty$, which completes the proof of \eqref{eq-pre-decay:Main}.
	\end{proof}
	
	From this proof, one can see that the derivatives of $\phi^+$ has actually exponential decay at infinity, while the polynomial decay is enough for the arguments in this article.

	\subsection{Conservation laws}\label{subsect-pre:Conservation}
	
	In this paragraph, we will rigorously prove the conservation of mass and energy introduced in Section~\ref{subsect-intro:ConservationHamilton}. These conservation laws are well-known for graph case where $\Sigma(t)=\{y=\eta(t,x)\}$, while the proof of general cases requires some geometrical arguments. To begin with, we introduce a technical lemma which will be frequently used in the proof of conservation laws and the virial identity \eqref{eq-intro-viri:Main} to be studied in Section~\ref{Sect:viri}.
	
	\begin{lemma}\label{lem-pre-cons:GenDeriFormula}
		Consider two domains $\Omega^\pm(t)$ filling the domain $\T^d\times]-H^-,H^+[$ and varying in time $t\in[0,T[$ with $H^\pm<+\infty$ and $T>0$. These domains are assumed to be separated by an interface $\Sigma(t)$ that admits a parametrization $\gamma:[0,T[\times\T^d\to \T^d\times]-H^-,H^+[$. We denote by $\tilde{\Omega}^\pm$ the time-space domains $\tilde{\Omega}^\pm = \{(t,x,y)\in [0,T[\times\mathbb{T}^d\times\mathbb{R} : (x,y)\in \Omega^\pm(t) \}$.
		
		If $\gamma$ is $C^1$ in space and time, then, for all functions $f^\pm \in C^1(\tilde{\Omega}^\pm)$, there holds
		\begin{equation}\label{eq-pre-cons:GenDeriFormula}
			\frac{d}{dt} \iint_{\Omega^\pm(t)} f^\pm dxdy = \iint_{\Omega^\pm(t)} \partial_t f^\pm dxdy \pm \int_{\Sigma(t)} n\cdot\gamma_t f^\pm dS.
		\end{equation}
		Note that $dS$ is the surface element on $\Sigma(t)$.
	\end{lemma}
	\begin{proof}
		To prove the formula \eqref{eq-pre-cons:GenDeriFormula}, we will prove that both side are equal as distributions in time, namely $\mathcal{D}'(]0,T[)$. This allows us to work on time-space domains, which is fixed. Then the desired equality follows from divergence theorem.
		
		Consider the time-space interface $\tilde{\Sigma}^\pm = \{(t,x,y)\in [0,T[\times\mathbb{T}^d\times\mathbb{R} : (x,y)\in \Sigma(t) \}$, which can be parameterized as $(t,x) \mapsto \left( t, \gamma(t,x) \right)$. Therefore, the tangent space of $\tilde{\Sigma}$ is spanned by $(1,\gamma_t)$ and $(0,\tau_j)$, where $(\tau_j)$ is arbitrary basis of the tangent space of $\Sigma(t)$. Consequently, the unit normal vector and surface element of $\tilde{\Sigma}$ read
		\begin{equation}\label{eq-pre-cons:TimeSpaceNor&Surf}
			\tilde{n} = \frac{(-n\cdot\gamma_t, n)}{\sqrt{1+(n\cdot\gamma_t)^2}},\ \ d\tilde{S} = \sqrt{1+(n\cdot\gamma_t)^2} dSdt,
		\end{equation}
		respectively. Here $n = n(t)$ denotes the unit normal vector of $\Sigma(t)$.
		
		For any test function $\varphi = \varphi(t)\in C^\infty_c(]0,T[)$, we have
		\begin{equation}\label{eq-pre-cons:GenDeriFormula-S1}
			\begin{aligned}
				\int_0^T \left( \frac{d}{dt} \iint_{\Omega^\pm(t)} f^\pm dxdy \right)& \varphi dt = - \iiint_{\tilde{\Omega}^\pm} f^\pm \partial_t\varphi dtdxdy \\
				=& \iiint_{\tilde{\Omega}^\pm} \partial_t f^\pm \varphi dtdxdy - \iiint_{\tilde{\Omega}^\pm} \partial_t (f^\pm \varphi) dtdxdy \\
				=& \int_0^T \iint_{\Omega^\pm(t)} \partial_t f^\pm \varphi dxdy dt - \iiint_{\tilde{\Omega}^\pm} \partial_t (f^\pm \varphi) dtdxdy.
			\end{aligned}
		\end{equation}
		To handle the second term on the right-hand side, we apply the divergence theorem in time-space,
		\begin{equation}\label{eq-pre-cons:GenDeriFormula-S2}
			\iiint_{\tilde{\Omega}^\pm} \partial_t (f^\pm \varphi) dtdxdy = \pm \iint_{\tilde{\Sigma}} \tilde{n} \cdot (f^\pm \varphi,0,0) d\tilde{S}.
		\end{equation}
		Note that the solid boundaries $\{y=\mp H^\pm\}$ does not appear since the time-space version of its normal vector equals $(0,\mp 1,0)$, whose scalar product with $(f \varphi,0,0)$ vanishes. By applying the formula \eqref{eq-pre-cons:TimeSpaceNor&Surf} of $\tilde{n}$ and $d\tilde{S}$, we have
		\begin{equation}\label{eq-pre-cons:GenDeriFormula-S3}
			\pm \iint_{\tilde{\Sigma}} \tilde{n} \cdot (f^\pm \varphi,0,0) d\tilde{S} = \mp \int_0^T\int_{\Sigma(t)} n\cdot\gamma_t f^\pm \varphi dS dt.
		\end{equation}
		By combining \eqref{eq-pre-cons:GenDeriFormula-S1}, \eqref{eq-pre-cons:GenDeriFormula-S2}, and \eqref{eq-pre-cons:GenDeriFormula-S3}, one concludes
		\begin{equation*}
			\int_0^T \left( \frac{d}{dt} \iint_{\Omega^\pm(t)} f^\pm dxdy \right) \varphi dt = \int_0^T \left( \iint_{\Omega^\pm(t)} \partial_t f^\pm dxdy \pm \int_{\Sigma(t)} n\cdot\gamma_t f^\pm dS \right) \varphi dt
		\end{equation*}
		for arbitrary $\varphi = \varphi(t)\in C^\infty_c(]0,T[)$, which implies the desired formula \eqref{eq-pre-cons:GenDeriFormula}.
	\end{proof}
	
	Now, we are in a position to prove the conservation laws.
	
	\begin{proposition}\label{prop-pre-cons:Main}
		Let $(\gamma, \phi^\pm,P^\pm)$ be a regular solution to the system \eqref{eq-intro-ww:Main} with depth $H^\pm\in]0,+\infty]$. Then the following conservation laws hold true:
		\begin{equation}\label{eq-pre-cons:Main}
			\frac{d}{dt} M = 0,\quad \frac{d}{dt} E = 0,
		\end{equation}
		where the mass $M$ and the total energy $E$ are defined in \eqref{eq-intro-ch:Mass} and \eqref{eq-intro-ch:TotEnergy}, respectively.
	\end{proposition}
	
	In order to overcome the difficulties caused by infinite depths, we will work on truncated domains defined as below.
	\begin{equation}\label{eq-pre-cons:ModDomain}
		\begin{aligned}
			&\Omega^+_{h} := \Omega^+ \cap \{y>-h^+\}, \quad \Omega^-_{h} := \Omega^- \cap \{y<h^-\}; \\
			&h^\pm = \left\{
			\begin{array}{ll}
				H^\pm, & \text{if }H^\pm<+\infty; \\
				\text{arbitrary number such that }\Sigma\subset\{\mp y <h^\pm\}, & \text{if }H^\pm=+\infty.
			\end{array}
			\right.
		\end{aligned}
	\end{equation}
	In the proof of the infinite depth case, we will work on these bounded domains and take the limit $h^\pm\to+\infty$ in the end. The proof of Proposition~\ref{prop-pre-cons:Main} will be divided into four steps: \textbf{(1)} prove conservation of mass; \textbf{(2)} compute the time derivative of potential energy; \textbf{(3)} compute the time derivative of truncated kinetic energy ; \textbf{(4)} conclude conservation of energy. 
	
	\medskip
	
	\noindent\textbf{Step 1: Conservation of mass.} Before calculating the time derivative of mass $M$ defined in \eqref{eq-intro-ch:Mass}, we give an alternative formulation of $M$. An application of divergence theorem in the domains $\Omega^\pm_{h}$ gives
	\begin{align*}
		\iint_{\Omega^\pm_{h}} dxdy = \iint_{\Omega^\pm_{h}} \nabla_{x,y}\cdot(0,y) dxdy =& \pm \int_\Sigma n\cdot(0,y) dS \mp \int_{y=\mp h^\pm} y dx \\
		=& \pm \int_\Sigma n\cdot(0,y) dS + (2\pi)^d h^\pm = \pm M + (2\pi)^d h^\pm,
	\end{align*}
	where the last equality follows from the definition \eqref{eq-intro-ch:Mass} of $M$. Thus, we conclude that
	\begin{equation}\label{eq-pre-cons:MassAlt}
		M = \pm \left( \iint_{\Omega^\pm_{h}} dxdy - (2\pi)^d h^\pm \right).
	\end{equation}
	
	From this formula, the time derivative of $M$ can be obtained via Lemma \ref{lem-pre-cons:GenDeriFormula} with $f=1$ and $\Omega^\pm$ replaced by $\Omega^\pm_{h}$.
	\begin{equation*}
		\frac{d}{dt} M = \pm \frac{d}{dt}\iint_{\Omega^\pm_{h}} dxdy = \int_{\Sigma} n\cdot\gamma_t dS.
	\end{equation*}
	Due to the third equation of \eqref{eq-intro-ww:Main}, $n\cdot\gamma_t$ can be replaced by $n\cdot\nabla_{x,y}\phi^\pm$. Then we apply the divergence theorem and use $\Delta_{x,y}\phi^\pm =0$ to compute
	\begin{equation*}
		\frac{d}{dt} M = \int_{\Sigma} n\cdot\nabla_{x,y}\phi^\pm dS = \int_{y=\mp h^\pm} \partial_y\phi^\pm dx.
	\end{equation*}
	In finite depth case, by choosing $h^\pm=H^\pm$ and using the boundary condition (the fifth equation in \eqref{eq-intro-ww:Main}), we conclude that the right-hand side is zero. In infinite depth case, as $h^\pm\to+\infty$, the the right-hand side tends to zero due to the rapid decay of $\partial_y\phi^\pm$ (see \eqref{eq-intro-instab:DecayScalPot}), which completes the proof of conservation law \eqref{eq-pre-cons:Main} for the mass $M$. \\

	\medskip
	
	\noindent\textbf{Step 2: Time derivative of potential energy.} From the definition \eqref{eq-intro-ch:PotEnergy} of the potential energy $E_p$, it contains two components concerning gravity and surface tension, respectively.
	\begin{align*}
		E_p =& E_p^\mathrm{grav} + E_p^\mathrm{surf}, \\
		E_p^\mathrm{grav} =& \frac{Ag}{2}  \int_{\Sigma} n\cdot(0,y^2) dS, \quad E_p^\mathrm{surf} = \frac{\sigma}{\rho^++\rho^-} \left[\mathrm{Area}(\Sigma(t)) - (2\pi)^d\right].
	\end{align*}
	In what follows, we will compute the time derivative of these two components respectively.
	
	\begin{lemma}\label{lem-pre-cons:TimeDeriPotEnergy}
		Under the hypothesis of Proposition~\ref{prop-pre-cons:Main}, the following formulas hold.
		\begin{align}
			\frac{d}{dt} E_p^\mathrm{grav} =& Ag \int_\Sigma n\cdot\gamma_t y dS, \label{eq-pre-cons:TimeDeriPotEnergyGrav} \\
			\frac{d}{dt} E_p^\mathrm{surf} =& \frac{\sigma}{\rho^++\rho^-} \int_{\Sigma} n\cdot\gamma_t \kappa dS. \label{eq-pre-cons:TimeDeriPotEnergySurf}
		\end{align}
		Recall that $\kappa$ is the mean curvature of the interface $\Sigma$.
	\end{lemma}
	
	\begin{proof}
		The time derivative \eqref{eq-pre-cons:TimeDeriPotEnergySurf} of $E_p^\mathrm{surf}$ can be simply obtained by the first variation of area formula (see, for instance, Theorem 5.19 of \cite{gallot1990riemannian}),
		\begin{equation*}
			\frac{d}{dt} \mathrm{Area}(\Sigma(t)) = \int_{\Sigma} n\cdot\gamma_t \kappa dS.
		\end{equation*}
		
		To calculate the time derivative of $E_p^\mathrm{grav}$, as in the proof of the conservation of mass, we could rewrite it as the linear combination of integrals on $\Omega^\pm_h$ (see \eqref{eq-pre-cons:ModDomain}) so that the formula \eqref{eq-pre-cons:GenDeriFormula} can be applied. The divergence theorem in $\Omega^\pm_{h}$ yields
		\begin{align*}
			2 \iint_{\Omega^\pm_{h}} y dxdy =& \iint_{\Omega^\pm_{h}} \nabla_{x,y}\cdot(0,y^2) dxdy \\
			=& \pm\int_{\Sigma} n\cdot(0,y^2) dS \mp \int_{y=\mp h^\pm} y^2 dx = \pm\int_{\Sigma} n\cdot(0,y^2) dS \mp (2\pi)^d(h^\pm)^2
		\end{align*}
		which implies
		\begin{equation}\label{eq-pre-cons:ReformGravPotEnergy}
			E_p^\mathrm{grav} = \sum_{\pm} \underline{\rho}^\pm g \left[ \iint_{\Omega^\pm_{h}} y dxdy \pm \frac{1}{2}(2\pi)^d(h^\pm)^2 \right].
		\end{equation}
		Recall that $\underline{\rho}^\pm = \rho^\pm/(\rho^++\rho^-)$. Then, via Lemma~\ref{lem-pre-cons:GenDeriFormula}, we can conclude \eqref{eq-pre-cons:TimeDeriPotEnergyGrav} by
		\begin{equation*}
			\frac{d}{dt} E_p^\mathrm{grav} = \sum_\pm \pm \underline{\rho}^\pm g \int_\Sigma n\cdot\gamma_t y dS = Ag \int_\Sigma n\cdot\gamma_t y dS,
		\end{equation*}
		where we used the fact that $A = \underline{\rho}^+ - \underline{\rho}^-$ (see \eqref{eq-intro-ch:AtwoodNum}).
	\end{proof}
	
	\begin{remark}\label{rmk-pre-cons:SignPotEnergy}
		Note that the contribution $E_p^\mathrm{grav}$ of gravity in potential energy has the same signature as the Atwood number $A$. To see this, we consider the following real numbers
		\begin{equation*}
			K^\pm := \iint_{\Omega^\pm_{h}} y dxdy \pm \frac{1}{2}(2\pi)^d(h^\pm)^2.
		\end{equation*}
		It is clear that $K^++K^- = 0$. Besides, we have
		\begin{equation*}
			K^+ \ge \iint_{\Omega^+_{h^+}\cap\{y<0\}} y dxdy + \frac{1}{2}(2\pi)^d(h^+)^2 \ge \iint_{\{y<0\}} y dxdy + \frac{1}{2}(2\pi)^d(h^+)^2 = 0.
		\end{equation*}
		These yield that there exists $K\ge 0$ such that $K^\pm = \pm K$. Thus, $E_p^\mathrm{grav}$ equals $AgK$, whose sign is the same as $A$. In particular, when $A\le 0$, $E_p^\mathrm{grav}$ becomes non-positive.
	\end{remark}

	\medskip
	
	\noindent\textbf{Step 3: Time derivative of truncated kinetic energy.} In order to calculate the time derivative of kinetic energy $E_k$ (see \eqref{eq-intro-ch:KinEnergy}), one observes that it is already written in terms of integrals on $\Omega^\pm$. Thus we can expect to apply the formula \eqref{eq-pre-cons:GenDeriFormula}, which is not valid for unbounded domains (namely, infinite depth case). To overcome this difficulty, we introduce a truncated kinetic energy,
	\begin{equation}\label{eq-pre-cons:KinEenergyCutOff}
		E_k^h := \frac{1}{2} \sum_{\pm} \underline{\rho}^\pm \iint_{\Omega^\pm_{h}} |\nabla_{x,y}\phi^\pm|^2 dxdy,
	\end{equation}
	where $\Omega^\pm_{h}$ are bounded domains defined by \eqref{eq-pre-cons:ModDomain}. Then we can compute the time derivative of $E_k^h$ by the following lemma.
	\begin{lemma}\label{lem-pre-cons:TimeDeriKinEnergy}
		Under the hypothesis of Proposition~\ref{prop-pre-cons:Main}, the following formula holds.
		\begin{equation}\label{eq-pre-cons:TimeDeriKinEnergy}
			\frac{d}{dt} E^h_k = - \frac{\sigma}{\rho^++\rho^-} \int_\Sigma n\cdot\gamma_t \kappa dS -Ag \int_\Sigma n\cdot\gamma_t y dS + \sum_\pm \mp \underline{\rho}^\pm \int_{y=\mp h^\pm} \phi^\pm_y \phi^\pm_t dx.
		\end{equation}
	\end{lemma}
	
	\begin{proof}
		Through Lemma~\ref{lem-pre-cons:GenDeriFormula}, the time derivative of truncated kinetic energy $E_k^h$, defined by \eqref{eq-pre-cons:KinEenergyCutOff}, reads
		\begin{equation*}
			\frac{d}{dt} E^h_k = \sum_\pm \underline{\rho}^\pm \left[ \iint_{\Omega^\pm_{h}} \nabla_{x,y}\phi^\pm \cdot \nabla_{x,y} \partial_t\phi^\pm dxdy \pm \int_{\Sigma}n\cdot\gamma_t \frac{|\nabla_{x,y}\phi^\pm|^2}{2} dS \right].
		\end{equation*}
		Then we apply divergence theorem to the integral in $\Omega^\pm_{h}$ together with $\Delta_{x,y}\phi^\pm=0$ (the second equation in \eqref{eq-intro-ww:Main}) and obtain
		\begin{align*}
			\iint_{\Omega^\pm_{h}} \nabla_{x,y}\phi^\pm \cdot \nabla_{x,y} \partial_t\phi^\pm dxdy =& \pm\int_\Sigma n\cdot\nabla_{x,y}\phi^\pm \partial_t\phi^\pm dS \mp \int_{y=\mp h^\pm} \phi^\pm_y \phi^\pm_t dx \\
			=& \pm\int_\Sigma n\cdot\gamma_t \partial_t\phi^\pm dS \mp \int_{y=\mp h^\pm} \phi^\pm_y \phi^\pm_t dx,
		\end{align*}
		where the second equality follows from the third equation in \eqref{eq-intro-ww:Main}. By combing the computations above, we can conclude by using the first equation in \eqref{eq-intro-ww:Main} that
		\begin{align*}
			\frac{d}{dt} E^h_k =& \sum_\pm \underline{\rho}^\pm \left[ \pm \int_{\Sigma}n\cdot\gamma_t \left( \phi^\pm_t + \frac{|\nabla_{x,y}\phi^\pm|^2}{2}\right) dS \mp \int_{y=\mp h^\pm} \phi^\pm_y \phi^\pm_t dx \right] \\
			=& \sum_\pm \left[ \mp \int_{\Sigma}n\cdot\gamma_t \left( \frac{P^\pm}{\rho^++\rho^-} + \underline{\rho}^\pm gy \right) dS \mp \underline{\rho}^\pm \int_{y=\mp h^\pm} \phi^\pm_y \phi^\pm_t dx \right] \\
			=&- \frac{\sigma}{\rho^++\rho^-} \int_\Sigma n\cdot\gamma_t \kappa dS -Ag \int_\Sigma n\cdot\gamma_t y dS + \sum_\pm \mp \underline{\rho}^\pm \int_{y=\mp h^\pm} \phi^\pm_y \phi^\pm_t dx.
		\end{align*}
		Here, we use the fourth equation in \eqref{eq-intro-ww:Main} to obtain the last equality.
	\end{proof}
	
	\medskip
	
	\noindent\textbf{Step 4: Conservation of energy.} Taking the sum of \eqref{eq-pre-cons:TimeDeriPotEnergyGrav}, \eqref{eq-pre-cons:TimeDeriPotEnergySurf}, and \eqref{eq-pre-cons:TimeDeriKinEnergy}, we have, for all time $t\in[0,T[$,
	\begin{equation*}
		E_p(t) + E_k^h(t) = E_p(0) + E_k^h(0) + \sum_\pm \mp \underline{\rho}^\pm \int_0^t\int_{\T^d} \phi^\pm_y(t',x,\mp h^\pm) \phi^\pm_t(t',x,\mp h^\pm) dx dt'.
	\end{equation*}
	Recall that, in finite depth case $H^\pm<+\infty$, we choose $h^\pm=H^\pm$ (see \eqref{eq-pre-cons:ModDomain}). Thus, the truncated kinetic energy $E_k^h$ equals the real one $E_k$ and the last integral on the right-hand side vanishes due to the fifth equation in \eqref{eq-intro-ww:Main}. The conservation of energy then follows.
	
	The main difficulty is the infinite depth case $H^+=+\infty$ or $H^-=+\infty$. Notice that, from the fifth equation of \eqref{eq-intro-ww:Main} and the condition (2) in Definition \ref{def-intro-instab:RegSol}, one could see that, as $h^\pm\to H^\pm$, the truncated kinetic energy $E^h_k$ tends to the real one $E_k$. Thus, it remains to show that, uniformly in time $t'\in[0,t]$, the integral
	\begin{equation*}
		\int_{\T^d} \phi^\pm_y(t',x,\mp h^\pm) \phi^\pm_t(t',x,\mp h^\pm) dx
	\end{equation*}
	tends to zero as $h^\pm\to+\infty$. Recall that the derivative $\phi^\pm_y$ of scalar potential has arbitrary polynomial decay at infinity \eqref{eq-intro-instab:DecayScalPot} (see also Lemma~\ref{lem-pre-decay:Main}), which reduces our problem to the uniform-in-time polynomial growth of $\phi^\pm_t$. From the first equation in \eqref{eq-intro-ww:Main}, we can rewrite $\phi^\pm_t$ as
	\begin{equation*}
		\phi^\pm_t = -\frac{|\nabla_{x,y}\phi^\pm|^2}{2} - \frac{P^\pm}{\rho^++\rho^-} - gy,
	\end{equation*}
	where the first term has polynomial decay thanks to \eqref{eq-intro-instab:DecayScalPot}, the last term is polynomial, and the growth of $P^\pm$ is assumed in the condition (6) in Definition~\ref{def-intro-instab:RegSol}. In conclusion, $\phi^\pm_y(\cdot,\mp h^\pm)\phi^\pm_t(\cdot,\mp h^\pm)$ decays polynomially at infinity, uniformly in time $t'\in[0,t]$, and its integral in $x$ tends to zero as $h^\pm\to+\infty$, which proves the conservation of energy in infinite depth case.

	\section{Virial theorem}\label{Sect:viri}
	
	The goal of this section is to complete the proof of Theorem~\ref{thm-intro-viri:Virial} via a direct calculus. Recall that the virial identity reads
	\begin{equation}\label{eq-viri:Main}
		\frac{1}{2}\frac{d}{dt} \int_{\Sigma}n\cdot(0,y\psi) dS = \tilde{E}_k - E_p + R, \quad R= R_b^+ + R_b^- + R_s,
	\end{equation}
	where $\psi$ is defined in \eqref{eq-intro-ch:DefPsi}, the modified kinetic energy $\tilde{E}_k$ is defined by \eqref{eq-intro-viri:ModKinEnergy} and the formula of reminders $R_b^\pm$ and $R_s$ is given by \eqref{eq-intro-viri:RemBottom} and \eqref{eq-intro-viri:RemSurfaceTension}, respectively. The main idea is similar to the proof of conservation laws in Section \ref{subsect-pre:Conservation}. We will rewrite the left-hand side of \eqref{eq-viri:Main} in terms of integrals on $\Omega^\pm$ (see Lemma~\ref{lem-viri:Reform}) and apply the formula \eqref{eq-pre-cons:GenDeriFormula}. The only difficulty is the term concerning curvature in \eqref{eq-viri:MainTemp}, which will be studied in Lemma~\ref{lem-viri:Curv}. To conclude Theorem~\ref{thm-intro-viri:Virial}, we also need to check that the reminders $R^\pm_b$ and $R_s$ are non-negative. The fact $R^\pm_b\ge 0$ is obvious, while $R_s \ge 0$ is not trivial and we shall use the same method as in Lemma~\ref{lem-viri:Curv} to prove it.
	
	As in Section~\ref{subsect-pre:Conservation}, to cover the infinite depth case $H^+=+\infty$ or $H^-=+\infty$, we will work on the truncated domains $\Omega^\pm_{h}$ defined by \eqref{eq-pre-cons:ModDomain}, which we recall below
	\begin{equation}\label{eq-viri:ModDomain}
		\begin{aligned}
			&\Omega^+_{h} := \Omega^+ \cap \{y>-h^+\}, \quad \Omega^-_{h} := \Omega^- \cap \{y<h^-\}; \\
			&h^\pm = \left\{
			\begin{array}{ll}
				H^\pm, & \text{if }H^\pm<+\infty; \\
				\text{arbitrary number such that }\Sigma\subset\{\mp y <h^\pm\}, & \text{if }H^\pm=+\infty,
			\end{array}
			\right.
		\end{aligned}
	\end{equation}
	and take the limit $h^\pm\to+\infty$ in the end.
	
	\begin{lemma}\label{lem-viri:Reform}
		Let $(\gamma,\phi^\pm,P^\pm)$ be a regular solution to \eqref{eq-intro-ww:Main} (see Definition~\ref{def-intro-instab:RegSol}) and $\Omega^\pm_{h}$ be truncated domains defined by \eqref{eq-viri:ModDomain}. Then there holds
		\begin{equation}\label{eq-viri:Reform}
			\int_{\Sigma}n\cdot(0,y\psi) dS = \sum_\pm \underline{\rho}^\pm \left[ \iint_{\Omega^\pm_{h}} \partial_{y}(y\phi^\pm) dxdy - h^\pm \int_{y=\mp h^\pm} \phi^\pm dx \right].
		\end{equation}
	\end{lemma}
	\begin{proof}
		Recall that $\psi$ is defined by \eqref{eq-intro-ch:DefPsi}. By divergence theorem on $\Omega^\pm_{h}$, we have
		\begin{align*}
			\int_{\Sigma}n\cdot(0,y\psi) dS =& \sum_\pm \underline{\rho}^\pm \int_{\Sigma}n\cdot(0,y\phi^\pm) dS \\
			=& \sum_\pm \underline{\rho}^\pm \left[ \iint_{\Omega^\pm_{h}} \nabla_{x,y}\cdot(0,y\phi^\pm) dxdy \pm \int_{y=\mp h^\pm} y\phi^\pm dx \right],
		\end{align*}
		which completes the proof.
	\end{proof}
	
	Based on Lemma~\ref{lem-viri:Reform}, it suffices to compute the time derivative of the terms on the right-hand side of \eqref{eq-viri:Reform}.
	\begin{lemma}\label{lem-viri:TimeDeri}
		Let $(\gamma,\phi^\pm,P^\pm)$ be a regular solution to \eqref{eq-intro-ww:Main} (see Definition~\ref{def-intro-instab:RegSol}) and $\Omega^\pm_{h}$ be the truncated domains defined by \eqref{eq-viri:ModDomain}. Then there hold
		\begin{align}
			&\begin{aligned}
				\frac{d}{dt} \iint_{\Omega^\pm_{h}} \partial_{y}(y\phi^\pm) dxdy = &\iint_{\Omega^\pm_{h}} \left( \frac{1}{2}|\nabla_x\phi^\pm|^2 + \frac{3}{2}|\partial_y\phi^\pm|^2 \right)dxdy \\
				&- \int_\Sigma n\cdot(0,y)\left(\pm\frac{P^\pm}{\rho^\pm}\right)dS - h^\pm\int_{y=\mp h^\pm} \frac{P^\pm}{\rho^\pm} dx \\
				& - 2g\iint_{\Omega^\pm_{h}}ydxdy \pm \int_{y=\mp h^\pm} \partial_y\phi^\pm \partial_{y}(y\phi^\pm) dx,
			\end{aligned} \label{eq-viri:TimeDeriMain} \\
			&\frac{d}{dt} \int_{y=\mp h^\pm} \phi^\pm dx = - \int_{y=\mp h^\pm} \frac{|\nabla_{x,y}\phi^\pm|^2}{2} dx - \int_{y=\mp h^\pm} \frac{P^\pm}{\rho^\pm} dx \pm gh^\pm (2\pi)^d. \label{eq-viri:TimeDeriBdy}
		\end{align}
	\end{lemma}
	\begin{proof}
		The second equality \eqref{eq-viri:TimeDeriBdy} can be obtained directly from the first equation of \eqref{eq-intro-ww:Main}. In what follows, we focus on the proof of the first equality \eqref{eq-viri:TimeDeriMain}. An application of the formula \eqref{eq-pre-cons:GenDeriFormula} gives
		\begin{align*}
			\frac{d}{dt} \iint_{\Omega^\pm_{h}} \partial_{y}(y\phi^\pm) dxdy =& \iint_{\Omega^\pm_{h}} \partial_{y}(y\phi^\pm_t) dxdy \pm \int_\Sigma n\cdot\gamma_t \partial_{y}(y\phi^\pm) dS \\
			=& - \iint_{\Omega^\pm_{h}} \partial_{y}\left(y \frac{|\nabla_{x,y}\phi^\pm|^2}{2}\right) dxdy - \iint_{\Omega^\pm_{h}} \partial_{y}\left(y \frac{P^\pm}{\rho^\pm}\right) dxdy \\
			&- \iint_{\Omega^\pm_{h}} \partial_{y}\left(gy^2\right) dxdy \pm \int_\Sigma n\cdot\nabla_{x,y}\phi^\pm \partial_{y}(y\phi^\pm) dS,
		\end{align*}
		where the second equality follows from the system \eqref{eq-intro-ww:Main}. Now, we apply divergence theorem to the second and the last term on the right-hand side,
		\begin{align*}
			\iint_{\Omega^\pm_{h}} \partial_{y}\left(y \frac{P^\pm}{\rho^\pm}\right) dxdy =& \pm\int_\Sigma n\cdot(0,y)\frac{P^\pm}{\rho^\pm}dS \mp \int_{y=\mp h^\pm} y\frac{P^\pm}{\rho^\pm} dx \\
			=& \int_\Sigma n\cdot(0,y)\left(\pm\frac{P^\pm}{\rho^\pm}\right)dS + h^\pm \int_{y=\mp h^\pm} \frac{P^\pm}{\rho^\pm} dx, \\
			\pm \int_\Sigma n\cdot\nabla_{x,y}\phi^\pm \partial_{y}(y\phi^\pm) dS  =& \iint_{\Omega^\pm_{h}} \nabla_{x,y}\cdot\left( \nabla_{x,y}\phi^\pm \partial_{y}(y\phi^\pm) \right) dxdy \\
			&\hspace{8em}\pm \int_{y=\mp h^\pm} \partial_y \phi^\pm \partial_{y}(y\phi^\pm) dx,
		\end{align*}
		and deduce that
		\begin{align*}
			\frac{d}{dt} \iint_{\Omega^\pm_{h}} \partial_{y}(y\phi^\pm) dxdy =& \iint_{\Omega^\pm_{h}} \left[ \nabla_{x,y}\cdot\left( \nabla_{x,y}\phi^\pm \partial_{y}(y\phi^\pm) \right) - \partial_{y}\left(y \frac{|\nabla_{x,y}\phi^\pm|^2}{2}\right) \right] dxdy \\
			& - \int_\Sigma n\cdot(0,y)\left(\pm\frac{P^\pm}{\rho^\pm}\right)dS - h^\pm \int_{y=\mp h^\pm} \frac{P^\pm}{\rho^\pm} dx \\
			& - 2g\iint_{\Omega^\pm_{h}} ydxdy\pm \int_{y=\mp h^\pm} \partial_{y}\phi^\pm \partial_{y}(y\phi^\pm) dx.
		\end{align*}
		Then the desired equality \eqref{eq-viri:MainTrivial} follows from the following identity,
		\begin{equation*}
			\nabla_{x,y}\cdot\left( \nabla_{x,y}\phi^\pm \partial_{y}(y\phi^\pm) \right) - \partial_{y}\left(y \frac{|\nabla_{x,y}\phi^\pm|^2}{2}\right) = \frac{1}{2}|\nabla_x\phi^\pm|^2 + \frac{3}{2}|\partial_{y}\phi^\pm|^2,
		\end{equation*}
		which can be computed directly by using $\Delta_{x,y}\phi^\pm=0$ (the second equation in \eqref{eq-intro-ww:Main}).
	\end{proof}
	
	By combining Lemma~\ref{lem-viri:Reform} and~\ref{lem-viri:TimeDeri}, we conclude that
	\begin{equation}\label{eq-viri:MainTrivial}
		\begin{aligned}
			\frac{d}{dt} \int_{\Sigma}n\cdot(0,y\psi) dS =& 2(\tilde{E}_k^h - E_p^\mathrm{grav} + R_{b,+}^h + R_{b,-}^h) \\
			& + \sum_{\pm} \pm \underline{\rho}^\pm \int_{y=\mp h^\pm} \partial_{y}\phi^\pm \partial_{y}(y\phi^\pm) dx - \frac{\sigma}{\rho^++\rho^-} \int_\Sigma n\cdot(0,y)\kappa dS,
		\end{aligned}
	\end{equation}
	where we used the fourth equation of \eqref{eq-intro-ww:Main} to make appear the mean curvature $\kappa$ and used \eqref{eq-pre-cons:ReformGravPotEnergy} to obtain $E_p^\mathrm{grav}$. Here $\tilde{E}_k^h$ and $R_{b,\pm}^h$ represent the truncated quantities
	\begin{equation*}
		\tilde{E}_k^h = \sum_{\pm} \underline{\rho}^\pm \iint_{\Omega^\pm_{h}} \left(\frac{1}{4}|\nabla_x\phi^\pm|^2 + \frac{3}{4}|\partial_y\phi^\pm|^2\right) dxdy, \quad R_{b,\pm}^h = \frac{1}{4} \underline{\rho}^\pm h^\pm \int_{y=\mp h^\pm}|\nabla_x\phi^\pm|^2dx,
	\end{equation*}
	whose limit as $h^\pm\to H^\pm$ equals $\tilde{E}_k$ (see \eqref{eq-intro-viri:ModKinEnergy}) and $R^\pm_b$ (see \eqref{eq-intro-viri:RemBottom}), respectively, according to the fifth equation in \eqref{eq-intro-ww:Main} and the decay \eqref{eq-intro-instab:DecayScalPot} of scalar potentials. With the same reason, the second term on the right-hand side of \eqref{eq-viri:MainTrivial} tends to zero as $h^\pm\to H^\pm$. As a conclusion, we have
	\begin{equation}\label{eq-viri:MainTemp}
		\frac{1}{2}\frac{d}{dt} \int_{\Sigma}n\cdot(0,y\psi) dS = \tilde{E}_k - E_p^\mathrm{grav} + R_b^+ + R_b^- - \frac{\sigma}{2(\rho^++\rho^-)} \int_\Sigma n\cdot(0,y)\kappa dS.
	\end{equation}
	In comparison with the target \eqref{eq-viri:Main}, it remains to compute the last integral on the right-hand side of \eqref{eq-viri:MainTemp}.
	\begin{lemma}\label{lem-viri:Curv}
		Let $(\gamma,\phi^\pm,P^\pm)$ be a regular solution to \eqref{eq-intro-ww:Main} (see Definition~\ref{def-intro-instab:RegSol}). Then there holds
		\begin{equation}\label{eq-viri:Curv}
			\int_\Sigma n\cdot(0,y)\kappa dS = \mathrm{Area}(\Sigma) - \int_\Sigma |n\cdot(0,1)|^2 dS.
		\end{equation}
	\end{lemma}
	
	To prove the identity \eqref{eq-viri:Curv}, we introduce another method to represent the fluid domains $\Omega^\pm$. For simplicity, we omit the dependence in time. There exists a $C^2$ function $\chi$ defined on $\Omega:=\T^d\times]-H^+,H^-[$, such that
	\begin{equation*}
		\Omega^+ = \{\chi>0\},\quad \Omega^-=\{\chi<0\}, \quad \Sigma = \{\chi=0\}.
	\end{equation*}
	Such $\chi$ is not unique, while our argument works for any choice of $\chi$. Let $F$ be an arbitrary continuous function defined on $\Omega$. Then the integral of $f:= F|_\Sigma$ on the interface $\Sigma$ equals
	\begin{equation*}
		\int_\Sigma f dS = \iint_{\Omega} F \delta(\chi) |\nabla_{x,y}\chi| dxdy := \lim_{\epsilon\to0} \iint_{\Omega} F \varphi_\epsilon(\chi) |\nabla_{x,y}\chi| dxdy,
	\end{equation*}
	where $\varphi_\epsilon := \epsilon^{-1}\varphi(\cdot/\epsilon)$ with $\varphi\in C^\infty_c(\R)$ and $\varphi(0)=1$ (cf. Theorem 6.1.5 of \cite{hormander1983analysis1}). Note that the normal direction $n$ and the mean curvature $\kappa$, which are both defined on $\Sigma$, could be rewritten as restrictions
	\begin{equation*}
		n = \left. -\frac{\nabla_{x,y}\chi}{|\nabla_{x,y}\chi|}\right|_{\Sigma}\text{ and } \kappa = \left. -\nabla_{x,y}\cdot\left(\frac{\nabla_{x,y}\chi}{|\nabla_{x,y}\chi|}\right)\right|_{\Sigma},
	\end{equation*}
	respectively.

	\begin{proof}[Proof of Lemma~\ref{lem-viri:Curv}]
		Through the representation above, the left-hand side of \eqref{eq-viri:Curv} reads
		\begin{equation}\label{eq-viri:Curv-S1}
			\begin{aligned}
				\int_\Sigma n\cdot(0,y)\kappa dS =& \iint_{\Omega} \frac{\nabla_{x,y}\chi}{|\nabla_{x,y}\chi|} \cdot (0,y) \nabla_{x,y}\cdot\left( \frac{\nabla_{x,y}\chi}{|\nabla_{x,y}\chi|} \right) \delta(\chi) |\nabla_{x,y}\chi| dxdy \\
				=& \iint_{\Omega} y\partial_{y}\chi \nabla_{x,y}\cdot\left( \frac{\nabla_{x,y}\chi}{|\nabla_{x,y}\chi|} \right) \delta(\chi) dxdy \\
				=& \lim_{\epsilon\to0} \iint_{\Omega} y\partial_{y}\chi \nabla_{x,y}\cdot\left( \frac{\nabla_{x,y}\chi}{|\nabla_{x,y}\chi|} \right) \varphi_\epsilon(\chi) dxdy.
			\end{aligned}
		\end{equation}
		Then we apply divergence theorem to the right-hand side (note that the integrands vanish near $y= \mp H^\pm$ due to the cut-off $\varphi_\epsilon$) and deduce that
		\begin{equation}\label{eq-viri:Curv-S2}
			\begin{aligned}
				&\hspace{-2em}\iint_{\Omega} y\partial_{y}\chi \nabla_{x,y}\cdot\left( \frac{\nabla_{x,y}\chi}{|\nabla_{x,y}\chi|} \right) \varphi_\epsilon(\chi) dxdy \\
				=& - \iint_{\Omega} \nabla_{x,y}(y\partial_{y}\chi)\cdot  \frac{\nabla_{x,y}\chi}{|\nabla_{x,y}\chi|}  \varphi_\epsilon(\chi) dxdy - \iint_{\Omega} y\partial_{y}\chi \frac{\nabla_{x,y}\chi}{|\nabla_{x,y}\chi|} \cdot \left(\varphi_\epsilon'(\chi) \nabla_{x,y}\chi\right) dxdy \\
				=& - \iint_{\Omega} \frac{(\partial_{y}\chi)^2}{|\nabla_{x,y}\chi|} \varphi_\epsilon(\chi) dxdy - \iint_{\Omega} \partial_y\left(|\nabla_{x,y}\chi|\right) y\varphi_\epsilon(\chi) dxdy  \\
				&\hspace{8em}- \iint_{\Omega} y\partial_{y}\chi|\nabla_{x,y}\chi| \varphi_\epsilon'(\chi) dxdy.
			\end{aligned}
		\end{equation}
		Via integration by parts, the second term on the right-hand side equals
		\begin{equation}\label{eq-viri:Curv-S3}
			- \iint_{\Omega} \partial_y\left(|\nabla_{x,y}\chi|\right) y\varphi_\epsilon(\chi) dxdy = \iint_{\Omega} |\nabla_{x,y}\chi| \varphi_\epsilon(\chi) dxdy + \iint_{\Omega} y\partial_{y}\chi|\nabla_{x,y}\chi| \varphi'_\epsilon(\chi) dxdy.
		\end{equation}
		
		By combing \eqref{eq-viri:Curv-S1}, \eqref{eq-viri:Curv-S2}, and \eqref{eq-viri:Curv-S3} above, we have
		\begin{align*}
			\int_\Sigma n\cdot(0,y)\kappa dS =& \lim_{\epsilon\to0} \left[- \iint_{\Omega} \frac{(\partial_{y}\chi)^2}{|\nabla_{x,y}\chi|} \varphi_\epsilon(\chi) dxdy + \iint_{\Omega} |\nabla_{x,y}\chi| \varphi_\epsilon(\chi) dxdy\right] \\
			=& -\iint_{\Omega} \left|\frac{\nabla_{x,y}\chi}{|\nabla_{x,y}\chi|} \cdot(0,1)\right|^2 \delta(\chi)|\nabla_{x,y}\chi| dxdy + \iint_{\Omega} \delta(\chi) |\nabla_{x,y}\chi| dxdy,
		\end{align*}
		where, in term of integrals on $\Sigma$, the right-hand side can be written as
		\begin{equation*}
			-\int_\Sigma |n\cdot(0,1)|^2 dS + \int_\Sigma dS,
		\end{equation*}
		which completes the proof of \eqref{eq-viri:Curv}.
	\end{proof}
	
	\begin{proof}[Proof of Theorem~\ref{thm-intro-viri:Virial}]
		By inserting \eqref{eq-viri:Curv} into \eqref{eq-viri:MainTemp} and using the definition \eqref{eq-intro-ch:PotEnergy} of potential energy $E_p$ together with the definition \eqref{eq-intro-viri:RemSurfaceTension} of $R_s$, we are able to conclude the virial identity \eqref{eq-viri:Main}. 
		\begin{align*}
			&\hspace{-2em}\frac{1}{2}\frac{d}{dt} \int_{\Sigma}n\cdot(0,y\psi) dS \\
			=& \tilde{E}_k - E_p^\mathrm{grav} + R_b^+ + R_b^- - \frac{\sigma}{2(\rho^++\rho^-)}\left[ \mathrm{Area}(\Sigma(t)) - \int_\Sigma |n\cdot(0,1)|^2 dS \right] \\
			=& \tilde{E}_k - E_p^\mathrm{grav} - E_p^\mathrm{surf} + R_b^+ + R_b^- \\
			&\hspace{4em}+ \frac{\sigma}{2(\rho^++\rho^-)}\left[ \mathrm{Area}(\Sigma(t)) + \int_\Sigma |n\cdot(0,1)|^2 dS - 2(2\pi)^d \right] \\
			=& \tilde{E}_k - E_p + R_b^+ + R_b^- + R_s.
		\end{align*}
		
		The only remaining point is to show that $R_s \ge 0$, which can be proved by using again the method in the proof of Lemma~\ref{lem-viri:Curv}. Recall that the remainder $R_s$ equals
		\begin{equation*}
			R_s = \frac{\sigma}{2(\rho^++\rho^-)} \left[ \int_\Sigma \left( 1 + |n\cdot(0,1)|^2 \right) dS - 2(2\pi)^d \right],
		\end{equation*}
		where the integral on $\Sigma$ can be expressed as an integral on $\Omega = \T^d\times]-H^+,H^-[$. Namely,
		\begin{align*}
			\int_\Sigma \left( 1 + |n\cdot(0,1)|^2 \right) dS =& \iint_{\Omega} \left( 1 + \left|\frac{\nabla_{x,y}\chi}{|\nabla_{x,y}\chi|}\cdot(0,1)\right|^2 \right) \delta(\chi)|\nabla_{x,y}\chi| dxdy \\
			=& \iint_{\Omega} \left( |\nabla_{x,y}\chi| + \frac{(\partial_{y}\chi)^2}{|\nabla_{x,y}\chi|} \right) \delta(\chi) dxdy.
		\end{align*} 
		Noticing that
		\begin{equation*}
			|\nabla_{x,y}\chi| + \frac{(\partial_{y}\chi)^2}{|\nabla_{x,y}\chi|} \ge 2|\partial_{y}\chi|,
		\end{equation*}
		we have
		\begin{align*}
			\int_\Sigma \left( 1 + |n\cdot(0,1)|^2 \right) dS \ge& 2\iint_{\Omega} |\partial_{y}\chi| \delta(\chi) dxdy \\
			=& 2\iint_{\Omega} \left|\frac{\nabla_{x,y}\chi}{|\nabla_{x,y}\chi|}\cdot(0,1)\right| \delta(\chi) |\nabla_{x,y}\chi| dxdy.
		\end{align*}
		In terms of integrals on $\Sigma$, this inequality can be rewritten as
		\begin{equation*}
			\int_\Sigma \left( 1 + |n\cdot(0,1)|^2 \right) dS \ge 2\int_\Sigma |n\cdot(0,1)| dS \ge 2\left|\int_\Sigma n\cdot(0,1) dS\right|.
		\end{equation*}
		Since $\nabla_{x,y}\cdot(0,1)=0$, through divergence theorem on $\Omega^+_h$ (see \eqref{eq-viri:ModDomain}), the last absolute value equals
		\begin{equation*}
			2\left|\int_\Sigma n\cdot(0,1) dS\right| = 2 \left| \int_{y=-h^+} dx \right| = 2(2\pi)^d,
		\end{equation*}
		which proves the desired result $R_s\ge 0$.
	\end{proof}

	\section{Instability of system}\label{Sect:instab}
	
	The purpose of this section is to prove the instability of the system \eqref{eq-intro-ww:Main} in graph case (Theorem~\ref{thm-intro-instab:Graph}) and 2D overlapping case (Theorem~\ref{thm-intro-instab:NonGraph}). The idea of the proof is similar to that in \cite{alazard2023virial}. Based on the virial identity \eqref{eq-intro-viri:Main}, we first show that the right-hand side can be bounded from below by the absolute value $|E|$ of the total energy $E$ (see \eqref{eq-intro-ch:TotEnergy}), which is preserved in time. Then we apply the trace estimate to control the integral on the left-hand side of \eqref{eq-intro-viri:Main} by the geometrical quantities such as slope and curvature.
	
	\subsection{The consequence of the virial identity}\label{subsect-ins:Viri}
	
	In this part, we aim to prove the following corollary of the virial identity \eqref{eq-intro-viri:Main}, stating that, in the setting of RTI or KHI problem, the integral on the left-hand side of \eqref{eq-intro-viri:Main} grows at least like $O(t)$.
	
	\begin{corollary}\label{cor-ins-viri:Main}
		Consider $d\ge 1$, $H^\pm\in]0,+\infty]$, $\rho^+\le\rho^-$, $\sigma\ge0$, and $T\in]0,+\infty]$. Let $(\gamma, \phi^\pm,P^\pm)$ be any regular solution to \eqref{eq-intro-ww:Main} (see Definition~\ref{def-intro-instab:RegSol}) on $[0,T[$. If the condition \eqref{eq-intro-instab:CondTotEnergy} recalled below holds
		\begin{equation*}
			E\neq 0\text{ when }\sigma=0,\text{ and }E<0\text{ when }\sigma>0,
		\end{equation*}
		then we have, for all time $t\in[0,T[$,
		\begin{equation}\label{eq-ins-viri:Main}
			I(t) \ge I(0) + |E|t,\quad I(t) = \int_{\Sigma}n\cdot(0,y\psi) dS,
		\end{equation}
		where $\psi$ is defined by \eqref{eq-intro-ch:DefPsi},
		$$
		\psi = \underline{\rho}^+\phi^+|_{\Sigma} - \underline{\rho}^-\phi^-|_{\Sigma},
		$$
		and $E$ is the total energy defined by \eqref{eq-intro-ch:TotEnergy}, which is preserved in time.
	\end{corollary}
	
	\begin{proof}
		By integrating in time on both sides of the virial identity \eqref{eq-intro-viri:Main}, we have
		\begin{equation*}
			I(t) = I(0) + 2\int_0^t \left( \tilde{E}_k(\tau) - E_p(\tau) + R(\tau) \right) d\tau \ge I(0) + 2\int_0^t \left( \frac{E_k(\tau)}{2} - E_p(\tau) \right) d\tau,
		\end{equation*}
		where the inequality follows from the fact that $R\ge 0$ and $\tilde{E_k} \ge E_k/2$, which can be seen by comparing the definition \eqref{eq-intro-ch:KinEnergy} and \eqref{eq-intro-viri:ModKinEnergy} of kinetic energy $E_k$ and modified kinetic energy $\tilde{E}_k$, respectively. Then the desired inequality \eqref{eq-ins-viri:Main} can be reduced to the following one,
		\begin{equation*}
			\frac{E_k(t)}{2} - E_p(t) \ge \frac{|E|}{2},\quad \forall t\in[0,T[
		\end{equation*}
		where $E$ is the total energy preserved in time. To prove this, we divide the condition \eqref{eq-intro-instab:CondTotEnergy} into two cases: (1) $E<0$; (2) $E>0$ with $\sigma=0$. For simplicity, we omit the dependence in time in the sequel.
		
		In the first case $E<0$, the proof is direct
		\begin{equation*}
			\frac{E_k}{2} - E_p = \frac{E_k}{2} - (E - E_k) = |E| + \frac{3}{2}E_k \ge |E|.
		\end{equation*}
		In the second case $E>0$ with $\sigma=0$, we observe from the definition \eqref{eq-intro-ch:PotEnergy} that the potential energy $E_p$ is non-negative since $\rho^+\le\rho^-$. Thus,
		\begin{equation*}
			\frac{E_k}{2} - E_p = \frac{1}{2}(E-E_p) - E_p = \frac{E}{2} - \frac{3}{2}E_p \ge \frac{E}{2} = \frac{|E|}{2},
		\end{equation*}
		which completes the proof.
	\end{proof}
	
	\subsection{The graph case}\label{subsect-ins:graph}
	
	In this section, we focus on the graph case, where the interface $\Sigma(t)$ is represented as $\{y=\eta(t,x)\}$ for some real-valued function $\eta$, and complete the proof of Theorem~\ref{thm-intro-instab:Graph}. By comparing the inequality \eqref{eq-ins-viri:Main} above and the target \eqref{eq-intro-instab:Graph}, it suffices to control the integral $I(t)$ defined in \eqref{eq-ins-viri:Main} above by a function of $\|\nabla_x\eta\|_{L^\infty}$, which will be achieved by trace estimate (see Proposition \ref{prop-ins-gr:Trace}). Since the time dependence has no impact in the arguments below, we will omit it for the sake of simplicity.
	
	To begin with, we reformulate $I(t)$ defined by \eqref{eq-ins-viri:Main} above in terms of $\eta$. By fixing the parametrization $\gamma$ as $\gamma(t,x) = (x,\eta(t,x))$, the normal direction $n$ and the surface element $dS$ read
	\begin{equation*}
		n = \frac{(-\nabla_x\eta, 1)}{\sqrt{1+|\nabla_x\eta|^2}},\quad dS = \sqrt{1+|\nabla_x\eta|^2} dx,
	\end{equation*}
	allowing us to express $I(t)$ as
	\begin{equation}\label{eq-ins-gr:ReformI}
		I = \int_{\T^d} \frac{(-\nabla_x\eta, 1)}{\sqrt{1+|\nabla_x\eta|^2}} \cdot (0,\eta\psi) \sqrt{1+|\nabla_x\eta|^2} dx = \int_{\T^d}\eta\psi dx.
	\end{equation}
	Recall that $\psi$ is defined by \eqref{eq-intro-ch:DefPsi}, which can be written as,
	\begin{equation}\label{eq-ins-gr:ReformPsi}
		\psi(x) = \underline{\rho}^+ \phi^+(x,\eta(x)) - \underline{\rho}^- \phi^-(x,\eta(x)).
	\end{equation}
	The main difficulty is to control $\psi$ by some quantities depending only on $\eta$. To do so, we introduce the following trace estimate, whose detailed proof can be found in Appendix B of \cite{alazard2023virial}.
	\begin{proposition}\label{prop-ins-gr:Trace}
		Let $(\gamma,\phi^\pm,P^\pm)$ be a regular solution to the system \eqref{eq-intro-ww:Main} (see Definition~\ref{def-intro-instab:RegSol}) with depths $H^\pm\in]0,+\infty]$. Then the following estimate holds,
		\begin{equation}\label{eq-ins-gr:Trace}
			\int_{\T^d} \left| \phi^\pm(x,\eta(x)) - m^\pm \right|^2 dx \le C (1+\|\nabla_x\eta\|_{L^\infty(\T^d)}) \int_{\Omega^\pm} \left| \nabla_{x,y}\phi^\pm \right|^2 dxdy,
		\end{equation}
		where $m^\pm\in\R$ is the average of $\phi^\pm(\cdot,\eta(\cdot))$ on $\T^d$, namely
		\begin{equation}\label{eq-ins-gr:DefMpm}
			m^\pm := \frac{1}{(2\pi)^d} \int_{\T^d} \phi^\pm(x,\eta(x)) dx.
		\end{equation}
		Here the constant $C>0$ depends on $d_0>0$ appearing in the condition (5) of Definition~\ref{def-intro-instab:RegSol}.
	\end{proposition}
	\begin{remark}
		Indeed, the left-hand side of \eqref{eq-ins-gr:Trace} can be refined as the Sobolev norm $H^{1/2}(\T^d)$, while $L^2(\T^d)$-norm is enough for the proof of Theorem~\ref{thm-intro-instab:Graph}.
	\end{remark}
	
	\begin{proof}[Proof of Theorem~\ref{thm-intro-instab:Graph}]
		By using the reformulation \eqref{eq-ins-gr:ReformI} of $I(t)$, the inequality \eqref{eq-ins-viri:Main} reads
		\begin{equation}\label{eq-ins-gr:MainPfS1}
			|E|t + \left.\int_{\T^d} \eta\psi dx\right|_{t=0} \le \int_{\T^d} \eta(t,x)\psi(t,x) dx.
		\end{equation}
		Recall that, in the condition (2) of Definition~\ref{def-intro-instab:RegSol}, we have assumed that $M=0$, where the mass $M$ is defined by \eqref{eq-intro-ch:Mass} and can be expressed, in the graph case, as
		\begin{equation}\label{eq-ins-gr:ReformZeroMass}
			M(t) = \int_{\T^d}\eta(t,x) dx \equiv 0.
		\end{equation}
		Therefore, from the formula \eqref{eq-ins-gr:ReformPsi} of $\psi$, the right-hand side of \eqref{eq-ins-gr:MainPfS1} equals
		\begin{equation*}
			\int_{\T^d} \eta(t,x)\psi(t,x) dx = \sum_{\pm} \pm \underline{\rho}^\pm \int_{\T^d} \eta(t,x)\left(\phi^\pm(t,x,\eta(x)) - m^\pm(t) \right) dx,
		\end{equation*}
		where the real number $m^\pm(t)$ (independent of $x$) is defined by \eqref{eq-ins-gr:DefMpm}. Thanks to the trace estimate \eqref{eq-ins-gr:Trace}, we can control the inner product of $\eta,\psi$ as
		\begin{equation}\label{eq-ins-gr:MainPfS2}
			\begin{aligned}
				\int_{\T^d}\eta\psi dx &\le \sum_{\pm} \underline{\rho}^\pm \|\eta\|_{L^2(\T^d)} \sqrt{\int_{\T^d} \left| \phi^\pm(x,\eta(x)) - m^\pm \right|^2 dx} \\
				&\lesssim_{d_0} \sum_{\pm} \|\eta\|_{L^2(\T^d)}\sqrt{ (1+\|\nabla_x\eta\|_{L^\infty(\T^d)}) \underline{\rho}^\pm \int_{\Omega^\pm} \left| \nabla_{x,y}\phi^\pm \right|^2 dxdy } \\
				&\lesssim \|\eta\|_{L^2(\T^d)} \sqrt{1+\|\nabla_x\eta\|_{L^\infty(\T^d)}} \sqrt{E_k},
			\end{aligned}
		\end{equation}
		where the last inequality follows from the definition \eqref{eq-intro-ch:KinEnergy} of kinetic energy $E_k$. Then, we replace the $\sqrt{E_k}$ on the right-hand side of \eqref{eq-ins-gr:MainPfS2} by $\sqrt{E-E_p}$. The potential energy $E_p$ defined by \eqref{eq-intro-ch:PotEnergy} can be expressed, in the graph case, as
		\begin{equation*}
			E_p = \frac{Ag}{2}\|\eta\|_{L^2(\T^d)}^2 + \frac{\sigma}{\rho^++\rho^-} \left[\mathrm{Area}(\Sigma(t)) - (2\pi)^d\right],
		\end{equation*}
		satisfying 
		\begin{equation*}
			E - E_p \le |E| - \frac{Ag}{2}\|\eta\|_{L^2(\T^d)}^2 \le |E| + |A|g\|\eta\|_{L^2(\T^d)}^2.
		\end{equation*}
		Note that the Atwood number $A = (\rho^+-\rho^-)/(\rho^++\rho^-)\in[-1,1]$ is negative when $\rho^+\le\rho^-$. Therefore, by combining \eqref{eq-ins-gr:MainPfS1} and \eqref{eq-ins-gr:MainPfS2}, we have
		\begin{equation}\label{eq-ins-gr:MainPfS3}
			|E|t + \left.\int_{\T^d} \eta\psi dx\right|_{t=0} \lesssim_{d_0} \|\eta\|_{L^2(\T^d)} \sqrt{1+\|\nabla_x\eta\|_{L^\infty(\T^d)}} \sqrt{|E|+|A|g\|\eta\|_{L^2(\T^d)}^2}.
		\end{equation}
		This inequality concludes the desired result \eqref{eq-intro-instab:Graph} by observing $\|\eta\|_{L^2(\T^d)} \lesssim \|\nabla_x\eta\|_{L^\infty(\T^d)}$, which is a consequence of the zero mass condition \eqref{eq-ins-gr:ReformZeroMass} (namely, the condition (3) in Definition~\ref{def-intro-instab:RegSol}).
	\end{proof}
	
	To end this part, let us assume that the regular solution $(\gamma,\phi^\pm,P^\pm)$ studied in Theorem~\ref{thm-intro-instab:Graph} is global-in-time (we do not know whether such solution exists). Directly from the inequality \eqref{eq-intro-instab:Graph}, we have 
	$$
	t^{2/5} \lesssim \|\nabla_x\eta(t)\|_{L^\infty(\T^d)}, \text{ as }t\to+\infty.
	$$
	In the absence of surface tension ($\sigma=0$), this growing rate of the slope $\|\nabla_x\eta(t)\|_{L^\infty(\T^d)}$ can be refined to $O(t^{2-\epsilon})$ for any $0<\epsilon\ll 1$. More precisely, we have
	\begin{proposition}\label{prop-ins-gr:GlobalGrow}
		Under the same assumption of Theorem~\ref{thm-intro-instab:Graph}, we further assume that $\sigma=0$ and the regular solution $(\gamma,\phi^\pm,P^\pm)$ is global-in-time. Then, for all $0<\epsilon\ll 1$, there holds
		\begin{equation}\label{eq-ins-gr:GlobalGrow}
			\limsup_{t\to+\infty} \frac{\|\nabla_x\eta(t)\|_{L^\infty(\T^d)}}{t^{2-\epsilon}} = +\infty.
		\end{equation}
	\end{proposition}
	\begin{proof}
		The proof is based on a refined version of Corollary~\ref{cor-ins-viri:Main}. We claim that, for all time $t\ge 0$, the right-hand side of the virial identity \eqref{eq-intro-viri:Main} can be bounded from below by
		\begin{equation}\label{eq-ins-gr:GlobalS1}
			\tilde{E}_k - E_p + R \ge \frac{|E| - E_p}{2}.
		\end{equation}
		To prove this, we first use the fact that $\tilde{E}_k\ge E_k/2$ and $R\ge 0$ and obtain
		$$ \tilde{E}_k - E_p + R \ge \frac{E_k}{2} - E_p. $$
		When $E>0$, the right-hand side reads
		$$ \frac{E}{2} - \frac{3}{2}E_p = \frac{|E|}{2} - \frac{3}{2}E_p \ge \frac{|E| - E_p}{2}, $$
		where we used the fact that $E_p\le 0$ since there in no surface tension $\sigma=0$. When $E<0$, we have
		$$ \frac{E}{2} - \frac{3}{2}E_p = \frac{E}{2} - (E-E_k) - \frac{E_p}{2} = \frac{|E|}{2} + E_k - \frac{E_p}{2} \ge \frac{|E| - E_p}{2}, $$
		which completes the proof of \eqref{eq-ins-gr:GlobalS1}. 
		
		By inserting \eqref{eq-ins-gr:GlobalS1} into the arguments in previous paragraphs, \eqref{eq-ins-gr:MainPfS3} becomes
		\begin{equation*}
			\int_0^t \frac{|E| - E_p(\tau)}{2} d\tau + I(0) \lesssim \|\eta\|_{L^2(\T^d)} \sqrt{1+\|\nabla_x\eta\|_{L^\infty(\T^d)}}\sqrt{E+|A|g\|\eta\|_{L^2(\T^d)}^2}.
		\end{equation*}
		Recall that $I(0)$, defined in \eqref{eq-ins-viri:Main}, depends only on initial data and $E_p = Ag\|\eta\|_{L^2(\T^d)}^2/2$ is non-positive. Therefore, the inequality above implies that 
		\begin{equation}\label{eq-ins-gr:GlobalS2}
			\int_0^t F(\tau)d\tau \le C \left( \sqrt{1+\|\nabla_x\eta(t)\|_{L^\infty(\T^d)}} F(t) + 1 \right),\quad F(t) = |E| + |A|g\|\eta(t)\|_{L^2(\T^d)}^2,
		\end{equation}
		for some the constant $C>0$ independent of time.
		
		Assume that \eqref{eq-ins-gr:GlobalGrow} is false. Namely, there exists time-independent constant $M>0$ such that $\|\nabla_x\eta(t)\|_{L^\infty(\T^d)}$ can be bounded by $Mt^{2-\epsilon}$ for all time $t\ge 1$. Then \eqref{eq-ins-gr:GlobalS2} implies
		\begin{equation}\label{eq-ins-gr:GlobalS3}
			\int_0^t F(\tau)d\tau \lesssim_M t^{1-\frac{\epsilon}{2}} \left( F(t) + 1 \right).
		\end{equation}
		We claim that this will lead to arbitrary polynomial growth of $F(t)$. Namely, for all $N\in\N$,
		\begin{equation}\label{eq-ins-gr:GlobalS4}
			t^{\frac{N\epsilon}{2}}\lesssim F(t),\quad \text{as }t\to+\infty,
		\end{equation}
		which is a contradiction since there exists time-independent constant $C_0>0$ such that, for all $N\in\N$,
		\begin{equation*}
			t^{\frac{N\epsilon}{2}}\lesssim F(t) \le |E| + |A|gC_0 \|\nabla_x\eta(t)\|_{L^\infty(\T^d)} \le |E| + |A|gC_0M t^{2-\epsilon},
		\end{equation*} 
		as $t\to+\infty$.
		
		It remains to prove \eqref{eq-ins-gr:GlobalS4} for all $N\in\N$ by induction. Firstly, the case $N=0$ is trivial from the definition \eqref{eq-ins-gr:GlobalS2} and the fact $E\neq 0$. Once \eqref{eq-ins-gr:GlobalS4} holds true for some $N\in\N$, the inequality \eqref{eq-ins-gr:GlobalS3} gives, 
		$$ t^{\frac{N\epsilon}{2}+1} \lesssim \int_0^t \tau^{\frac{N\epsilon}{2}}d\tau \lesssim \int_0^t F(\tau)d\tau \lesssim t^{1-\frac{\epsilon}{2}} (F(t)+1),  $$
		which yields
		$$ t^{\frac{(N+1)\epsilon}{2}} \lesssim F(t)+1,\quad\text{as }t\to+\infty. $$
		This is equivalent to \eqref{eq-ins-gr:GlobalS4} with $N$ replaced by $N+1$. Then the desired estimate \eqref{eq-ins-gr:GlobalS4} follows from an induction in $N$.
	\end{proof}

	\subsection{The 2D overlapping case}\label{subsect-ins:NonGraph}
	
	In this part, we shall focus on the 2D overlapping case and complete the proof of Theorem~\ref{thm-intro-instab:NonGraph}. The strategy of the proof is the same as that in Section~\ref{subsect-ins:graph}. The only difference is the trace estimate \eqref{eq-ins-gr:Trace}, while, in 2D overlapping case, we need to replace the slope $\|\nabla_x\eta\|_{L^\infty}$ appearing on the right-hand side of \eqref{eq-ins-gr:Trace} by some quantities depending on the curvature $\kappa$.
	
	In 2D overlapping case, the interface can be expressed via arc-length parameter $s$. In what follows, we fix the parametrization $\gamma(t)$ as
	\begin{equation}\label{eq-ins-ngr:ArcLength}
		\begin{array}{lccc}
			\gamma(t) \colon & \R/L(t) & \to & \T\times]-H^+,H^-[ \\[0.5ex]
			& s & \mapsto & \gamma(t,s) = (\alpha(t,s),\beta(t,s)),
		\end{array}
	\end{equation}
	where $L(t)$ is the length of the interface $\Sigma(t)$ and $\alpha,\beta$ are horizontal and vertical components of $\gamma(t)$, respectively. Recall that the interface is assumed to satisfy the \textit{chord-arc} condition (the condition (4) of Definition~\ref{def-intro-instab:RegSol}),
	\begin{equation}\label{eq-ngr:chord-arc}
		c_0|s-s'| \leqslant |\gamma(t,s)-\gamma(t,s')| \leqslant |s-s'|,\ \ \forall s,s'\in\mathbb{R}/L(t).
	\end{equation}
	Here $|s-s'|$ should be understood as the distance in $\mathbb{R}/L(t)$ and $c_0>0$ is a universal constant. 
	
	In terms of arc-length parametrization $\gamma(t,s)$, the normal direction reads $n = (-\beta_s,\alpha_s)$, and the curvature $\kappa$ is given by
	\begin{equation}\label{eq-ins-ngr:Curvature}
		\tau_s(s) = -\kappa(s)n(s),
	\end{equation}
	where $\tau := \gamma_s = (\alpha_s, \beta_s)$ is the unit tangent direction. Besides, the surface element $dS = ds$, yielding that the mass $M$ defined by \eqref{eq-intro-ch:Mass}, the potential energy $E_p$ defined by \eqref{eq-intro-ch:PotEnergy}, and the integral $I$ defined in \eqref{eq-ins-viri:Main} can be written as
	\begin{equation}\label{eq-ins-ngr:Reform}
		\begin{aligned}
			M =& \int_0^L \beta\alpha_s ds, \quad E_p = \frac{Ag}{2}\int_0^L \beta^2\alpha_sds + \frac{\sigma}{\rho^++\rho^-}(L-2\pi),\\
			I =& \int_0^L \beta \alpha_s \psi ds.
		\end{aligned}
	\end{equation}
	Recall that $\psi$ is defined by \eqref{eq-intro-ch:DefPsi} 
	$$
	\psi(t,s) = \underline{\rho}^+ \phi^+(t,\gamma(t,s)) - \underline{\rho}^- \phi^-(t,\gamma(t,s))
	$$
	in terms of the traces of scalar potentials $\phi^\pm$ at the interface $\Sigma$.
	
	Now we introduce the trace estimate in 2D overlapping case, which is an analogue of Proposition~\ref{prop-ins-gr:Trace}. Since this result is independent of time, we will omit the time dependence for simplicity.
	\begin{proposition}\label{prop-ins-ngr:Trace}
		Let $(\gamma,\phi^\pm,P^\pm)$ be a regular solution to the system \eqref{eq-intro-ww:Main} with depths $H^\pm\in]0,+\infty]$ (see Definition~\ref{def-intro-instab:RegSol}). We define $\epsilon>0$ as in \eqref{eq-intro-instab:DefEpsilon},
		\begin{equation}\label{eq-ins-ngr:DefEpsilon}
			\epsilon := \min \left(\frac{c_0}{N_0 (\|\kappa\|_{L^\infty}+1)},\ d_0 \right),
		\end{equation}
		where $0<c_0\le 1$ is the chord-arc constant defined in \eqref{eq-ngr:chord-arc}, $d_0>0$ is the constant appearing in the condition (5) of Definition~\ref{def-intro-instab:RegSol} (in infinite depth case, $d_0=+\infty$), and $N_0 \gg 1$ is a large enough universal constant. For all $C^1$ function $F$ on $\Omega^\pm$ with $\|\nabla_{x,y}F\|_{L^2(\Omega^\pm)}<+\infty$, if its trace $f = F|_\Sigma$ has zero mean
		\begin{equation*}
			\int_{\Sigma} f ds = 0,
		\end{equation*}
		then there exists a universal constant $C>0$ such that
		\begin{equation}\label{eq-ins-ngr:TraceGen}
			\|f\|_{L^2(\Sigma)} \le C \frac{L}{\sqrt{\epsilon}} \|\nabla_{x,y}F\|_{L^2(\Omega^\pm)}.
		\end{equation}
		Recall that $L$ is the length of the interface $\Sigma$. 
	\end{proposition}
	
	\begin{proof}
		We focus on the domain $\Omega^+$, while the other one $\Omega^-$ can be treated in the same way. The proof of this proposition is based on flattening by \textit{normal geodesic coordinate} defined by
		\begin{equation}\label{eq-ngr:nor-geo-coord}
			\begin{array}{lccc}
				\Phi \colon &\mathbb{R}/L \times [-\epsilon,\epsilon] & \rightarrow & \mathbb{T}\times]-H^+,H^-[ \\[0.5ex]
				&(s,r) & \mapsto & \gamma(s) + r n(s) 
			\end{array}
		\end{equation}
		It is clear that, for $N_0$ large enough, $\Phi$ is diffeomorphic (see Proposition~\ref{prop-nor:Main}) and its Jacobian is of determinant $j(s,r)=1+r\kappa(s) \in [1/2,3/2]$, since
		\begin{equation*}
			|r\kappa(s)| \leqslant \epsilon \|\kappa\|_{L^\infty} \leqslant \frac{c_0}{N_0} \leqslant \frac{1}{N_0} \ll 1.
		\end{equation*}
		Remark that $\{r<0\}\subset\Omega^+$, $\{r=0\} = \Sigma$, and $\{r>0\}\subset\Omega^-$. In the sequel, we shall use the restriction
		\begin{equation*}
			-\epsilon \leqslant r \leqslant 0.
		\end{equation*}
		
		To begin with, we truncate $F$ inside $\{-\epsilon<r\leqslant 0\}$ using
		\begin{equation*}
			G(s,r) := F\circ\Phi(s,r) \chi\left(\frac{r}{\epsilon}\right),
		\end{equation*}
		where $\chi\in C^\infty_{c}(\mathbb{R})$ is supported inside $]-1,1[$ and equals to $1$ near zero. One may observe that $G$ has the same trace as $F$ and 
		\begin{align*}
			\iint |\partial_r G(s,r)|^2 dsdr =& \iint \left| n(s)\cdot\left(\nabla_{x,y}F\right)\circ\Phi(s,r) \chi\left(\frac{r}{\epsilon}\right) + \frac{1}{\epsilon}F\circ\Phi(s,r)\chi'\left(\frac{r}{\epsilon}\right)  \right|^2 dsdr \\
			\lesssim& \iint_{\Omega^+} \left| \nabla_{x,y}F \right|^2 dxdy + \frac{1}{\epsilon^2} \int_{-\epsilon}^0\int_0^L \left| F\circ\Phi \right|^2 dsdr. 
		\end{align*}
		
		Since $G$ is supported in $\{-\epsilon < r \leqslant 0\}$, we have trivially,
		\begin{equation*}
			f(s) = G(s,0) = \int_{-\epsilon}^0 \partial_r G (s,r) dr,
		\end{equation*}
		thus
		\begin{equation*}
			|f(s)|^2 \leqslant \epsilon \int_{-\epsilon}^0 |\partial_r G(s,r)|^2 dr.
		\end{equation*}
		By integrating in $s$ on both sides, one obtains
		\begin{equation*}
			\|f\|_{L^2(\Sigma)}^2 \leqslant \epsilon \iint |\partial_r G|^2 dsdr \lesssim \epsilon\iint_{\Omega^+} \left| \nabla_{x,y}F \right|^2 dxdy + \frac{1}{\epsilon} \int_{-\epsilon}^0\int_0^L \left| F\circ\Phi \right|^2 dsdr.
		\end{equation*}
		
		In order to conclude \eqref{eq-ins-ngr:TraceGen}, we try to apply Poincar\'e's inequality for the last term on right-hand side. The difficulty is that the mean value $m$ of $F\circ\Phi$ in $[0,L]\times [-\epsilon,0]$, namely
		\begin{equation*}
			m = \frac{1}{\epsilon L}\int_{-\epsilon}^{0}\int_{0}^{L} F\circ\Phi(s,r) dsdr,
		\end{equation*}
		is non zero in general. To handle this problem, it suffices to check that this mean value $m$ can also be bounded by $\|\nabla_{x,y}F\|_{L^2(\Omega)}$. In fact, by repeating the argument above and applying the zero mean condition for $f=F\circ\Phi(\cdot,0)$, one may write $m$ as an integral of $\partial_r(F\circ\Phi)$,
		\begin{align*}
			m =& \frac{1}{\epsilon L}\int_{-\epsilon}^0\int_0^L F\circ\Phi(s,r) dsdr \\
			=& \frac{1}{\epsilon L}\int_{-\epsilon}^0\int_0^L \left( F\circ\Phi(s,0) + \int_{0}^{r} \partial_r (F\circ\Phi)(s,r') dr' \right) dsdr \\
			=& \frac{1}{\epsilon L}\int_{-\epsilon}^0\int_0^L \int_{0}^{r} \partial_r (F\circ\Phi)(s,r') dr'dsdr,
		\end{align*}
		which implies the desired control of $m$ :
		\begin{align*}
			|m| \leqslant& \frac{\sqrt{\epsilon}}{\epsilon L} \int_{-\epsilon}^0\int_0^L \sqrt{\int_{-\epsilon}^0 |\partial_r (F\circ\Phi)(s,r')|^2 dr'} dsdr \\
			\leqslant& \frac{ \sqrt{\epsilon} \sqrt{\epsilon L} }{\epsilon L} \sqrt{\int_{-\epsilon}^0\int_0^L\int_{-\epsilon}^0 |\partial_r (F\circ\Phi)(s,r')|^2 dr'dsdr} \\
			=& \sqrt{\frac{\epsilon}{L}}\|\partial_r (F\circ\Phi)\|_{L^2([0,L]\times[-\epsilon,0])}  \lesssim \sqrt{\frac{\epsilon}{L}} \|\nabla_{x,y} F\|_{L^2(\Omega^+)}
		\end{align*}
		
		Now we are able to apply Poincar\'e's inequality to $F\circ\Phi-m$ on $[0,L]\times[-\epsilon,0]$ : 
		\begin{align*}
			\int_{-\epsilon}^0\int_0^L \left| F\circ\Phi \right|^2 dsdr \leqslant& 2\int_{-\epsilon}^0\int_0^L \left| F\circ\Phi-m \right|^2 dsdr + 2\epsilon L m^2 \\
			\lesssim& \max\{L,\epsilon\}^2 \int_{-\epsilon}^\epsilon\int_0^L \left| \nabla_{s,r}(F\circ\Phi) \right|^2 dsdr + \epsilon^2 \|\nabla_{x,y} F\|_{L^2(\Omega^+)}^2 \\ \lesssim& L^2  \|\nabla_{x,y} F\|_{L^2(\Omega^+)}^2,
		\end{align*}
		where the last inequality follows from the fact that $L\geqslant 2\pi$ and $\epsilon \ll 1$ by choosing $N_0 \gg 1$.
	\end{proof}
	
	Recall that the aim of Theorem~\ref{thm-intro-instab:NonGraph} is to find a lower bound of $\epsilon$. Therefore, it is essential to control $L$, which appears on the right-hand side of \eqref{eq-ins-ngr:TraceGen}, in terms of $\epsilon$. In infinite depth case, this is impossible since it is easy to construct examples with large $L$ and small $\epsilon$, or with small $L$ and large $\epsilon$. However, under finite depth assumption $H^\pm<+\infty$, the following lemma allows us to control $L$ by $\epsilon^{-1}$.
	
	\begin{lemma}\label{lem-ngr:length-curvature}
		Under the assumptions of Proposition~\ref{prop-ins-ngr:Trace}, if we further assume that $H^\pm < +\infty$, then the following inequality holds true.
		\begin{equation}\label{eq-ngr:length-curvature}
			L\epsilon \leqslant 2\pi(H^++H^-).
		\end{equation}
	\end{lemma}
	\begin{proof}
		To prove this inequality, we shall use the change of variable $\Phi$ defined in \eqref{eq-ngr:nor-geo-coord},
		\begin{equation*}
			\Phi(s,r) = \gamma(s) + r n(s).
		\end{equation*}
		Let us denote by $\mathcal{O}_\epsilon$ the image of $[0,L]\times[-\epsilon,\epsilon]$ by $\Phi$, which is a subdomain of $\Omega^+\cup\Sigma\cup\Omega^-$. Therefore,
		\begin{equation*}
			2\pi(H^++H^-) = \text{Vol}(\Omega^+\cup\Sigma\cup\Omega^-) \geqslant \text{Vol}(\mathcal{O}_\epsilon) = \int_0^L\int_{-\epsilon}^\epsilon j(s,r) dsdr \geqslant L\epsilon,
		\end{equation*}
		where $j$ is the determinant of the Jacobian of $\Phi$ and, due to the definition \eqref{eq-ngr:nor-geo-coord} of $\epsilon$,
		$$
		j(s,r) = 1+ r\kappa(s) \ge 1- \epsilon \|\kappa\|_{L^\infty} \ge \frac{1}{2},
		$$
		when $N_0\gg1$ is chosen large enough.
	\end{proof}
	
	Now, we are at the position to prove the main result Theorem~\ref{thm-intro-instab:NonGraph}. As Theorem~\ref{thm-intro-instab:Graph}, we will combine Corollary~\ref{cor-ins-viri:Main}, the trace estimate \eqref{eq-ins-ngr:TraceGen}, and the inequality \eqref{eq-ngr:length-curvature}.
	
	\begin{proof}[Proof of Theorem~\ref{thm-intro-instab:NonGraph}]
		By comparing the consequence \eqref{eq-ins-viri:Main} of virial identity and the target inequality \eqref{eq-intro-instab:NonGraph}, it suffices to control the integral $I$ defined in \eqref{eq-ins-viri:Main} (see also its reformulation \eqref{eq-ins-ngr:Reform}) by the right-hand side of \eqref{eq-intro-instab:NonGraph}. Since the arguments below are independent of time, we shall omit the time dependence for simplicity.
		
		Recall that in the condition (3) of Definition~\ref{def-intro-instab:RegSol}, we have assumed that the total mass $M$ is identically zero. Thus, from the reformulation \eqref{eq-ins-ngr:Reform} and the definition \eqref{eq-intro-ch:DefPsi} of $\psi$, we can write the integral $I$ as
		\begin{equation*}
			I = \sum_\pm \pm \underline{\rho}^\pm \int_0^L \beta\alpha_s \phi^\pm|_{\Sigma} ds =  \sum_\pm \pm \underline{\rho}^\pm \int_0^L \beta\alpha_s (\phi^\pm|_{\Sigma} - m^\pm) ds,
		\end{equation*}
		where $m^\pm\in\R$ is the mean value of $\phi^\pm$ on $\Sigma$,
		\begin{equation*}
			m^\pm = \frac{1}{L}\int_0^L \phi^\pm(\gamma(s)) ds.
		\end{equation*}
		Clearly, $\phi^\pm-m^\pm$ has zero mean on the interface $\Sigma$, allowing us to apply the trace estimate \eqref{prop-ins-ngr:Trace}.
		\begin{align*}
			I &\le \sum_\pm \underline{\rho}^\pm \|\beta\alpha_s\|_{L^2(\Sigma)} \|\phi^\pm|_{\Sigma} - m^\pm\|_{L^2(\Sigma)} \\
			&\lesssim \sum_\pm \underline{\rho}^\pm \|\beta\alpha_s\|_{L^2(\Sigma)} L \epsilon^{-1/2} \|\nabla_{x,y}\phi^\pm\|_{L^2(\Omega^\pm)} \\
			&\lesssim L \epsilon^{-1/2} \|\beta\alpha_s\|_{L^2(\Sigma)} \sqrt{E_k},
		\end{align*}
		where the last inequality follows from the definition \eqref{eq-intro-ch:KinEnergy} and the fact $\underline{\rho}^\pm = \rho^\pm/(\rho^++\rho^-)\in[0,1]$. On the right-hand side, the $L^2(\Sigma)$-norm of $\beta\alpha_s$ can be estimated by observing $|\beta|\le L$ and $|\alpha_s|\le |\gamma_s| = 1$,
		\begin{equation*}
			\|\beta\alpha_s\|_{L^2(\Sigma)}^2 = \int_0^L |\beta|^2 |\alpha_s|^2 ds \le L^2 \int_0^L ds = L^3,
		\end{equation*}
		while $E_k$ can be estimated through the expression \eqref{eq-ins-ngr:Reform} of $E_p$,
		\begin{equation*}
			E_k = E - E_p \le E + \frac{|A|g}{2} \int_0^L \beta^2 \alpha_s ds \le |E| + |A|g L^3.
		\end{equation*}
		
		In conclusion, we have shown that the integral $I$ can be bounded by
		\begin{equation*}
			I \lesssim L^{5/2} \epsilon^{-1/2} \sqrt{|E| + |A|g L^3},
		\end{equation*}
		which, thanks to Lemma~\ref{lem-ngr:length-curvature}, implies
		\begin{equation*}
			I \le C \epsilon^{-3} \sqrt{|E| + |A|g \epsilon^{-3}},
		\end{equation*}
		for some constant $C>0$ depending on $H^\pm$. The proof of \eqref{eq-intro-instab:NonGraph} is completed.
	\end{proof}

	\appendix
	
	\section{Normal geodesic coordinate}\label{App:nor-geo}
	
	The purpose of this section is the following proposition, stating that the mapping $\Phi$, known as \textit{normal geodesic coordinate}, defined in \eqref{eq-ngr:nor-geo-coord}, realizes a diffeomorphism to its image.
	
	\begin{proposition}\label{prop-nor:Main}
		Consider $L>0$, $0<c_0\le 1$, and $N_0 >0$. Let $\gamma: \R/L \to \T\times\R$ be the arc-length parametrization of a $C^2$ curve $\Sigma$ with length $L$, which satisfies the chord-arc condition: for all $s_1,s_2\in \R/L$,
		\begin{equation}\label{eq-nor:ChordArc}
			c_0|s_1-s_2| \le |\gamma(s_1) - \gamma(s_2)| \le |s_1-s_2|.
		\end{equation}
		We define $\epsilon\in]0,+\infty]$ as 
		\begin{equation}\label{eq-nor:DefEspsilon}
			\epsilon := \frac{c_0}{N_0 \|\kappa\|_{L^\infty}},
		\end{equation}
		where $\kappa$ is the curvature of the curve $\Sigma$. Then, if $N_0$ is large enough, the mapping
		\begin{equation}\label{eq-nor:DefPhi}
			\begin{array}{lccc}
				\Phi : &\mathbb{R}/L \times ]-\epsilon,\epsilon[ & \rightarrow & \mathbb{T}\times\mathbb{R} \\[0.5ex]
				&(s,r) & \mapsto & \gamma(s) + rn(s)
			\end{array}
		\end{equation}
		is a diffeomorphism to its image, where $n(s)$ is the unit normal vector.
	\end{proposition}
	
	\begin{remark}
		Note that the choice of $\epsilon$ above is not the same as \eqref{eq-ins-ngr:DefEpsilon}. The extra $d_0$ in \eqref{eq-ins-ngr:DefEpsilon} is to make sure that the image of $\Phi$ lies in the fluid domain $\T\times]-H^+,H^-[$, especially in finite depth case. Actually, once Proposition~\ref{prop-nor:Main} is proved, the injectivity of $\Phi$ can be obtained for any choice of $\epsilon$ smaller than \eqref{eq-nor:DefEspsilon} above.
	\end{remark}
	
	If $\|\kappa\|_{L^\infty}=0$, the curve $\Sigma$ is straight. In this case we have $\epsilon=+\infty$ and that $\Phi$ is trivially a bijecttion. For non-trivial case $\|\kappa\|_{L^\infty}>0$ (thus $\epsilon\in]0,+\infty[$), one can compute the determinant $j$ of the Jacobian of $\Phi$,
	$$ j(s,r) = 1+r\kappa(s) \in [1/2,3/2], $$
	when $N_0$ is large enough. Thus, $\Phi$ is at least a local diffeomorphism with the desired correspondence. In the sequel, we shall check that this property is global by showing the injectivity of $\Phi$. If $\Phi$ is not injective, then
	\begin{equation}\label{eq-nor-geo:assum}
		\exists \text{ distinct points }(s,r),(s',r')\in [0,L]\times[-\epsilon,\epsilon],\text{ such that }\gamma(s) + rn(s) = \gamma(s') + r'n(s').
	\end{equation}
	Without loss of generality, we may assume $s<s'$, since $s=s'$ implies $r=r'$ by definition of $\Phi$. We will first prove that \eqref{eq-nor-geo:assum} holds only when $s$ and $s'$ are close enough (see Lemma~\ref{lem-nor-geo:local} below) and deduce Proposition~\ref{prop-nor:Main} by showing that $\Phi$ must be locally injective.
	
	In this section, the identification that $[0,2\pi]\times\mathbb{R}$ is a subdomain of $\mathbb{C}$ will also be used, which enables us to define a real function $\theta$ such that
	\begin{equation}\label{eq-nor:DefTheta}
		\gamma_s(s) = \left(\alpha_s(s), \beta_s(s)\right) = e^{i\theta(s)}.
	\end{equation}
	Then, by definition, the curvature reads
	\begin{equation}\label{eq-nor:ReformCurvature}
		\kappa = -\tau_s\cdot n = -\gamma_{ss} \cdot n = - \Imag(\gamma_{ss} \overline{i\gamma_s}) = -\theta'.
	\end{equation}
	
	\subsection{Reduction to local injectivity.}
	
	\begin{lemma}\label{lem-nor-geo:local}
		Let $N_0 \geqslant 2$ and $\|\kappa\|_{L^\infty}>0$. Then, under the assumption of Proposition~\ref{prop-nor:Main}, \eqref{eq-nor-geo:assum} implies
		\begin{equation}\label{eq-nor-geo:local}
			|s-s'|\leqslant \frac{4}{N_0 \|\kappa\|_{L^\infty}}.
		\end{equation}
	\end{lemma}
	\begin{proof}
		In fact, the assumption \eqref{eq-nor-geo:assum} can be written as
		\begin{equation*}
			\gamma(s')-\gamma(s) = r n(s) - r' n(s') = (r-r')n(s) + r'\int_{s'}^s \kappa(s_1)\tau(s_1) ds_1,
		\end{equation*}
		where $\tau$ is the unit tangent vector of $\Sigma$. The length of left-hand side can be bounded from below by applying chord-arc condition \eqref{eq-nor:ChordArc},
		\begin{equation*}
			c_0|s-s'| \leqslant |\gamma(s')-\gamma(s)|,
		\end{equation*}
		while the length of right-hand side is controlled by
		\begin{equation*}
			|r-r'| + \left| r'\int_{s'}^s \kappa(s_1)\tau(s_1) ds_1 \right| \leqslant |r-r'| + |s-s'|\epsilon\|\kappa\|_{L^\infty} \leqslant |r-r'| + \frac{c_0}{N_0}|s-s'|.
		\end{equation*}
		By combining the estimates above, we have
		\begin{equation*}
			c_0|s-s'| \leqslant |r-r'| + \frac{c_0}{N_0}|s-s'|.
		\end{equation*}
		Then, we can conclude by using the condition $N_0\geqslant 2$,
		\begin{equation*}
			|s-s'| \le \frac{2}{c_0}|r-r'| \le \frac{4}{c_0}\epsilon \le \frac{4}{N_0 \|\kappa\|_{L^\infty}},
		\end{equation*}
		which gives the desired inequality \eqref{eq-nor-geo:local}.
	\end{proof}
	
	\subsection{Local injectivity.} We shall check below that, when $s$ and $s'$ are close in the sense of \eqref{eq-nor-geo:local}, \eqref{eq-nor-geo:assum} can not hold, which is a contradiction and Proposition~\ref{prop-nor:Main} follows.
	\begin{lemma}\label{lem-nor:local}
		Let $N_0 \gg 1$ and $\|\kappa\|_{L^\infty}>0$. Then, under the assumption of Proposition~\ref{prop-nor:Main}, the inequality \eqref{eq-nor-geo:local} yields
		\begin{equation*}
			\tau(s)\cdot\tau(s') \geqslant \frac{1}{2}.
		\end{equation*}
	\end{lemma}
	\begin{proof}
		By identifying the unit tangent vector $\tau$ as complex number (see \eqref{eq-nor:DefTheta}), we have
		\begin{equation*}
			\tau(s)\cdot\tau(s') = \cos(\theta(s)-\theta(s')).
		\end{equation*}
		Then the expression \eqref{eq-nor:ReformCurvature} of the curvature $\kappa$ gives
		\begin{equation*}
			\left|\theta(s)-\theta(s')\right| \leqslant |s-s'| \max_{[s,s']}|\theta'| = |s-s'|\max_{[s,s']}|\kappa| \leqslant \frac{4}{N_0 \|\kappa\|_{L^\infty}} \|\kappa\|_{L^\infty} = \frac{4}{N_0}.
		\end{equation*}
		By choosing $N_0$ large enough, we can assume $\cos(4/N_0) \geqslant 1/2$, which proves the desired result.
	\end{proof}
	
	\begin{proof}[Proof of Proposition~\ref{prop-nor:Main}]
		Now we are able to deduce a contradiction to \eqref{eq-nor-geo:assum}, based on Lemma~\ref{lem-nor:local}. Let us consider the real-valued function : 
		\begin{equation*}
			\begin{array}{lccc}
				f : &[-\epsilon,\epsilon] & \rightarrow & \mathbb{R} \\[0.5ex]
				&h & \mapsto & \left( \gamma(s') + h n(s') - \gamma(s) \right)\cdot \tau(s)
			\end{array}.
		\end{equation*}
		Once \eqref{eq-nor-geo:assum} holds true, $f$ will admit a zero $h=r'$, while
		\begin{align*}
			f(h) =& \left( \gamma(s') + h n(s') - \gamma(s) - h n(s) \right)\cdot \tau(s) \\
			=& \left( \int_s^{s'}\tau(s_1) ds_1 + h \int_s^{s'}n_s(s_1) ds_1 \right)\cdot \tau(s) \\
			=& \left( \int_s^{s'}\tau(s_1) ds_1 + h \int_s^{s'}\kappa(s_1)\tau(s_1) ds_1 \right)\cdot \tau(s) \\
			=& \int_s^{s'} \left(1+ h\kappa(s_1)\right)\tau(s_1) \cdot \tau(s) ds_1.
		\end{align*}
		From Lemma~\ref{lem-nor:local}, $\tau(s_1) \cdot \tau(s) \geqslant 1/2$, and the definition of $\epsilon$ ensures that 
		$$1+ h\kappa(s_1) \geqslant 1-\epsilon\|\kappa\|_{L^\infty} \geqslant 1-\frac{c_0}{N_0} \geqslant 1/2,$$
		for large enough $N_0$. Therefore, $f$ has a lower bound $|s-s'|/4>0$, which is a contradiction.
	\end{proof}

	\printindex
	
	\printbibliography[heading=bibliography,title=References]

\end{document}